# On developing piecewise rational mapping with fine regulation capability for WENO schemes

Qin Li, Pan Yan, Xiao Huang, Liuqing Yang, Fengyuan Xu

School of Aerospace Engineering, Xiamen University, Xiamen, Fujian, 361102, China

**Abstract**: On the idea of mapped WENO-JS scheme, properties of mapping methods are analyzed, uncertainties in mapping development are investigated, and new rational mappings are proposed. Based on our former understandings, i.e. mapping at endpoints {0, 1} tending to identity mapping, an integrated $C_{m,n}$ condition is summarized for function development. Uncertainties, i.e., whether the mapping at endpoints would make mapped scheme behave like WENO or ENO, whether piecewise implementation would entail numerical instability, and whether WENO3-JS could preserve the third-order at first-order critical points by mapping, are analyzed and clarified. A new piecewise rational mapping with sufficient regulation capability is developed afterwards, where the flatness of mapping around the linear weights and its endpoint convergence toward identity mapping can be coordinated explicitly and simultaneously. Hence, the increase of resolution and preservation of stability can be balanced. Especially, concrete mappings are determined for WENO3,5,7-JS. Numerical cases are tested for the new mapped WENO-JS, which regards numerical stability including that in long time computation, resolution and robustness. In purpose of comparison, some recent mappings such as IM by [App. Math. Comput. 232, 2014:453-468], RM by [J. Sci. Comput. 67, 2016:540-580] and AIM by [J. Comput. Phys. 381, 2019:162-188] are chosen; in addition, some recent WENO-Z type scheme are selected also. Proposed new schemes can preserve optimal orders at corresponding critical points, achieve numerical stability and indicate overall comparative advantages regarding accuracy, resolution and robustness.

**Keywords**: WENO; Mapping method

## 1 Introduction

As one of the popular high-order difference methods, WENO schemes [1-2] especially WENO-JS in Ref. 2 have achieved large success. In order to improve performances, endeavors have been made constantly, which at least lie in two aspects, i.e. the mapping method [3] and WENO-Z approach [4-6]. In this study, further studies regarding the former method are concerned.

In Ref. 3, Henrick et al. examined the accuracy relation of WENO5-JS and proposed corresponding sufficient and necessary condition as well as the sufficient condition [3] to achieve the optimal fifth-order. They pointed out that when $f' = 0$ happened at $x_j$ or nearby, the above conditions would be violated and order degradation occurred. For remedy, they pioneered the idea of mapping method by introducing $g(\omega)$ for nonlinear weight $\omega \in [0,1]$, and new non-normalized weights were derived accordingly. A concrete *rational* function, referred as $g_M$ herein, was proposed and corresponding WENO5-M scheme was obtained. It is worthwhile to mention that Henrick et al. [3] and subsequent authors [7] did not think WENO3-JS could preserve its optimal order at the first-order critical point by mapping.

In Ref. 7, Feng et al. studied the stability and errors of WENO5-M in the case of long time computation. They found the scheme be liable to errors near discontinuities, and the reason was considered as the enlargement of nonlinear weights by $g_M$ near $\omega = 0,1$. For improvement, the authors suggested the mapped scheme to preform like ENO by appealing to the condition $g'(0,1) = 0$. A concrete piecewise *polynomial* mapping was proposed [7] as PM with the order $n+2$, where $n$ means $\text{PM}^{(i)}(d_k) = 0$ for $1 \leq i \leq n$ and $\text{PM}^{(n+1)}(d_k) \neq 0$. The choice of $n=6$ was suggested for WENO5 [7]. Still out of the similar concern, Feng et al. later devised a family of specific *rational* polynomial mapping called as IM [8] with free parameter $A$ and order degree $n$ just referred. It is noted that parameter recommendations of IM for WENO5 [8] would yield highly amplified weights near $\omega = 0,1$ than that in $g_M$, which is opposed in Ref. 7.

In Ref. 9, motivated by enhancing the resolution of WENO-JS by increasing the flatness of mapping around $\omega = d_k$, Li et al. made investigation on *polynomial* mapping independently. They increased the degree of polynomial degree to engender higher-order critical point at $\omega = d_k$ and flatter mapping profile there. It was found that [9] the choice pf single polynomial would be liable to ill-defined function with oscillations, whereas the piecewise polynomials (PPM) could easily fulfill the job. Although the error in long-time computation [7] was not referred there, the manner of polynomial approaching endpoints was concerned for the sake of numerical stability. Consequently, the approach to identity mapping at endpoints was chosen and the following conditions was proposed [9]: $g'(0,1) = 1$ and $g^{(i)}(0,1) = 0$ with $2 \leq i \leq m$ for certain $m$. The conditions indicate the mapped scheme resembles WENO near endpoints. As an example, a sixth-order PPM was proposed [9].

In order to realize the potential of *single rational* mapping functions with ENO-like feature at endpoints, Wang et al. [10] further studied a specific formulation called as RM with assigned orders of critical points at $\omega_k = 0$ and $d_k$. In order to obtain well-defined rational mapping, they claimed: (1) the condition such as $\{g^{(i)} = 0, 1 \leq i \leq m\}$ can only be imposed for one endpoint, where the point at $\omega = 0$ is chosen as the appropriate; (2) aforementioned $m$ and $n$ should be even where $n$ denotes the order of critical point at $\omega = d_k$. A choice of $(m, n) = (2, 6)$ is suggested. One can verify that the flatness of *RM* at $\omega = d_k$ and its convergence rate to ENO at $\omega = 0$ can only be adjusted through $m$ and $n$, and moreover the two features are competitive to each other. Recently, Vevek et al. [11] improved IM [8] by upgrading $A^{-1}\omega(1-\omega)$ in the mapping to $A^{-1}\left[\omega(1-\omega)\right]^{m-1}$, through which the new function can satisfy the conditions at endpoints proposed by Ref. 9. Furthermore, they proposed an adaptation algorithm as that in Ref. 12 to evaluate $A$, and the so-called AIM mapping for WENO7-JS was obtained. The adaptation of $A$ is to make mapping have flat profile at $\omega = d_k$ in smooth solution, and tend to identity mapping near endpoints at discontinuities.

In spite of progresses, uncertainties still exist. For example, (1) Feng et al. [7-8, 10] suggested ENO rather than WENO to tend to by mapping near endpoints, whereas Li [9] indicated the alternative manner of WENO there. Further investigations are preferred in this regard. (2) Wang et al. [10] deemed piecewise mapping have errors liable to numerical instability, and therefore the mapping function should be infinitely smooth or the single form is preferred. The authors even exemplified an occurrence of oscillations by WENO9 plus PM as the evidence. Whether piecewise mapping would be less appropriate than single one need further clarification. (3) Current studies except that of AIM [11] have indicated that the flatness of mapping at $\omega_k = d_k$ and convergence to identity mapping at endpoints are competitive or unachievable simultaneously. While for AIM, on the one hand the adaptive implementation indicates the absence of explicit and fine control, on the other hand its robustness and applicability for schemes such as WENO5-JS is unclear. It is obvious that other novel methods are expected. (4) Current studies considered the optimal order of WENO3 could not be preserved by mapping at critical points, however, which is rather needed by engineering. Targeting at above issues, clarification studies are made in this study and new piecewise rational mappings are proposed with explicit and fine control.

The paper is arranged as follows: typical mapping methods for WENO-JS are first reviewed in Section 2; uncertainties mentioned above are discussed in Section 3; in Section 4, new piecewise rational mappings are proposed and intensive analysis is made; careful numerical validations are carried out in Section 5; at last conclusions are drawn in Section 6.

## 2 Mapping function method for WENO-JS

In order to facilitate discussion, the complete set of WENO-JS with the orders from 3-9 is revisited first. Then, typical mapping methods are reviewed.

### 2.1 WENO-JS formula [2]

Consider the following one-dimensional hyperbolic conservation law

$$u_t + f(u)_x = 0. \tag{1}$$

Suppose the grids are equally partitioned as $x_j = j\Delta x$ where $\Delta x$ denotes the interval and $j$ is the grid index, Eq. (1) at $x_j$ can be re-written in conservative form as: $(u_t)_j = -(h_{j+1/2} - h_{j-1/2})/\Delta x$, where $f(x) = \frac{1}{\Delta x}\int_{x-\Delta x/2}^{x+\Delta x/2} h(x')dx'$. If $h(x)$ is approximated by $\hat{f}(x)$, the semi-discretized conservative scheme is written as

$$(u_t)_j = -(\hat{f}_{j+1/2} - \hat{f}_{j-1/2})/\Delta x. \tag{2}$$

To facilitate the description of WENO-JS formula, $f(u)_x > 0$ is assumed tentatively. Suppose $r$ is the number of substencils and also grid points of each one, WENO-JS [2] with the order $2r-1$ at smooth region can be formulated as:

$$\hat{f}_{j+1/2} = \sum_{k=0}^{r-1} \omega_k q_k^r \quad \text{with} \quad q_k^r = \sum_{l=0}^{r-1} a_{kl}^r f(u_{j-r+k+l+1}) \tag{3}$$

where $q_k^r$ is the candidate scheme with coefficients $a_{kl}^r$ shown in Table 8 in Appendix I and $\omega_k$ is the nonlinear weight. $\omega_k$ is derived from corresponding linear weight $d_k^r$ (shown in Table 9 in Appendix I) as

$$\omega_k = \alpha_k \Big/ \sum_{l=0}^{2} \alpha_l \quad \text{with} \quad \alpha_k = d_k^r \Big/ \left(\varepsilon + IS_k^{(r)}\right)^2 \tag{4}$$

where $\alpha_k$ denotes the non-normalized weight, $\varepsilon = 10^{-6}$ in WENO-JS and $IS_k^{(r)}$ is smoothness indicator. For WENO-JS, $IS_k^{(r)}$ in positive semi-definite quadratic form can be formulated as:

$$IS_k^{(r)} = \sum_{m=0}^{r-2} c_m^r \left( \sum_{l=0}^{r-1} b_{kml}^r f(u_{j-r+k+l+1}) \right)^2 \tag{5}$$

where the coefficients $b_{kml}^r$, $c_m^r$ are tabulated in Table 10-11 in Appendix I [2, 10-12].

2.2 Mapping method and functions

(1) Principle and properties of mapping

As stated in Ref. 3, when critical points occur in computation, the necessary and sufficient condition to acquire the order 2r-1 is hard to satisfy by WENO-JS, and accuracy loss ensues. However, the optimal order can still be achieved if the following sufficient condition be satisfied somehow

$$\omega_k = d_k + O(\Delta x^r), \tag{6}$$

where the superscript "r" of $d_k^r$ has been dropped for brevity. In order to achieve Eq. (6), a recipe was proposed by Henrick et al., in which prerequisite of accuracy relation between $\omega_k$ and $d_k$ should stand. The least requirement is:

$$\omega_k = d_k + O(\Delta x). \tag{7}$$

With Eq. (7), if one function $g(\omega_k)$ would make

$$g(\omega_k) = d_k + O(\omega_k - d_k)^r, \tag{8}$$

then $g(\omega_k)$ can serve as the non-normalized weight, and its normalized value can make Eq. (6) established [3]. By Taylor expansion, the requirement of $g(\omega_k)$ can be derived as [3]:

$$g(d_k) = d_k \text{ and } g^{(i)}(d_k) = 0 \text{ for } i < r. \tag{9}$$

In Ref. 3, Eq. (9) together with boundary conditions, namely $\{g(0,1) = 0,1\}$, were proposed as the requirements for mapping. However, in pursue of additional benefits as shown in introduction, more properties are preferred. Based on existing understandings, the properties concerned by mapping can be summarized as:

(a) Accuracy requirements such as Eq. (9) and boundary condition should be satisfied [3].
(b) Monotonicity, which makes the ordinal relation of original weights preserved.
(c) Flatness of mapping profile. Ref. 9 indicated that the increase of flatness at $\omega_k = d_k$

could enhance the resolution of WENO, which was realized by asking *higher* order of critical point than that in Eq. (8). Besides of the referred flatness, another kind is implied in previous investigations [8, 11] at the region away from $\omega_k = d_k$. As will be shown later (see Fig. 5), although a mapping function with higher order of critical point at $\omega_k = d_k$ would appear flatter there than another one with lower order, it is probable that the former might appear less flat when away from $\omega_k = d_k$ than the latter. For convenience, we refer the flatness in the neighborhood of $\omega_k = d_k$ as Flatness-I and the one away from the position as Flatness-II in this study.

(d) Convergence of mapping towards endpoints. As not involved in Ref. 3, the issue was first concerned by Refs. 7 and 9 separately. There are two different patterns of convergences: (i) the mapping converges to endpoints by $g^{(i)}(0,1) = 0$ with $1 \leq i \leq m$ and the mapped scheme is expected to behave like ENO there [7]; (ii) the mapping behaves by $\{g'(0,1) = 1, g^{(i)}(0,1) = 0 \text{ for } 1 < i \leq m\}$ and corresponding scheme tends to WENO-JS [8] at endpoints. Further discussions on two patterns will be given in Section 3, and for convenience the property of endpoint convergence is abbreviated as *PEC*. Ref. 10 indicated that it was difficult to acquire a realization of the first pattern with higher order *m*.

Based on above discussions, the conditions to define a mapping function $g(\omega)$ with given orders of derivatives at $\{0, d_k, 1\}$, which is referred as $C_{n,m,k}$, are proposed as a summary:

$$g^{(i)}(d_k) = \begin{cases} d_k, & i = 0 \\ 0, & 1 \leq i \leq n \\ \neq 0, & i = n+1 \end{cases}, \quad g^{(i)}(0) = \begin{cases} 0, & i = 0 \\ 1 | 0, & i = 1 \text{ if } m \geq 1 \\ 0, & 2 \leq i \leq m \\ \neq 0, & i = m+1 \end{cases}, \text{ and } g^{(i)}(1) = \begin{cases} 1, & i = 0 \\ 1 | 0, & i = 1 \text{ if } k \geq 1 \\ 0, & 2 \leq i \leq k \\ \neq 0, & i = k+1 \end{cases}. \quad (10)$$

When $m=k$, $C_{n,m,k}$ is simplified as $C_{n,m}$. Particularly, $g^{(1)}(0,1) = 1$ in Eq. (10) corresponds to the first pattern of convergence [7] above and $g^{(1)}(0,1) = 0$ corresponds the second [9]. It is conceivable that the larger *n*, *m* and *k* are, the higher rate of convergence two flatness and *PEC* will take.

(2) Mapping functions

In this point, typical mapping functions including our recent practices will be reviewed. To clearly demonstrate the relationship of functions with $C_{n,m}$ or $C_{n,m,k}$, their names will take the form as: $g_{n,m}^N$ or $g_{n,m,k}^N$, where *N* represents the order of polynomial or that of nominator in the rational situation. And for piecewise mapping, the definition takes: $g(\omega) = \{g^L, 0 < \omega < d_k; g^R, d_k \leq \omega < 1\}$. In the following, the subscript "*k*" of $\omega_k$ will be dropped for convenience. In addition, $\varepsilon$ in Eq. (4) will take $10^{-40}$ after taking mapping unless otherwise noted.

(a) $\text{PM}_{n,1}^{n+2}$ [7]

The mapping is a piecewise polynomial with the form

$$\begin{cases} (\text{PM}^L)_{n,1}^{n+2} = d_k + \frac{(-1)^n (n+1)}{d_k^{n+1}} (\omega - d_k)^{n+1} \left(\omega + \frac{d_k}{n+1}\right) \\ (\text{PM}^R)_{n,1}^{n+2} = d_k - \frac{n+1}{(1-d_k)^{n+1}} (\omega - d_k)^{n+1} \left(\omega + \frac{d_k - (n+2)}{n+1}\right) \end{cases},$$

which satisfies $C_{n,1}$ condition with $g'(0,1) = 0$. The concrete suggestion for WENO5-9 is $\text{PM}_{6,1}^{8}$ by Ref. 7, which was referred as PM6 therein. $\varepsilon$ in Eq. (4) takes $10^{-40}$ for mapped WENO-JS by the reference.

(b) $\text{PPM}_{n,m}^{n+m+1}$ [9]

The mapping is another kind of piecewise polynomial that is derived by employment of a general polynomial $\sum_{i=0}^{n+m+1} a_i \omega^i$ to satisfy $C_{n,m}$ in piecewise manner, where $g'(0,1) = 1$ in the case of $m \geq 1$. The specific cases of $m = 0$ at $n = 1\ldots 4$ and $m = 1$ at $n = 4$ were investigated in Ref. 9. One can see that $\text{PPM}_{n,1}^{n+2}$ would be analogous to $\text{PM}_{n,1}^{n+2}$ [7] in having the same orders in $C_{n,m}$ but differs in *PEC*. Based on studies in Ref. 7, we further propose the following general form of $(\text{PPM})_{n,m}^{n+m+1}$ as:

$$\begin{cases} (\text{PPM}^L)_{n,m}^{n+m+1} = d_k + \frac{(-1)^{n+m}}{d_k^{n+m}} (\omega - d_k)^{n+1} \sum_{i=0}^{m} a_i^m \omega^{m-i} d_k^i \\ (\text{PPM}^R)_{n,m}^{n+m+1} = d_k + \frac{1}{(1-d_k)^{n+m}} (\omega - d_k)^{n+1} \sum_{i=0}^{m} a_i^m (1-\omega)^{m-i} (1-d_k)^i \end{cases},$$

where $a_i^m = \frac{\prod_{j=0}^{m-1-i} (n+j)}{(m-i)!}$ for $i < m$ and $a_m^m = 1$. By means of the following theorem with its proof in Appendix III, $(\text{PPM}^R)_{n,m}^{n+m+1}$ satisfies $C_{n,m}$. Similarly, one can find that the counterpart $(\text{PPM}^L)_{n,m}^{n+m+1}$ in $[0, d_k]$ satisfies $C_{n,m}$ also.

**Theorem 1**: Considering a function as $f(\omega) = d_k + \frac{1}{(1-d_k)^{n+m}} (\omega - d_k)^{n+1} \sum_{i=0}^{m} a_i^m (1-\omega)^{m-i} (1-d_k)^i$ in $[d_k, 1]$ where $a_i^m = \frac{\prod_{j=0}^{m-1-i} (n+j)}{(m-i)!}$ for $i < m$ and $a_m^m = 1$, then $f(\omega)$ satisfies $C_{n,m}$.

To facilitate coding, the concrete forms of $(\text{PPM})_{n,m}^{n+m+1}$ with $m = 0\ldots 3$ are given as follows:

$$\begin{cases} (\text{PPM}^L)_{n,0}^{n+1} = d_k + \frac{(-1)^n}{d_k^n} (\omega - d_k)^{n+1} \\ (\text{PPM}^R)_{n,0}^{n+1} = d_k + \frac{1}{(1-d_k)^n} (\omega - d_k)^{n+1} \end{cases}, \quad \begin{cases} (\text{PPM}^R)_{n,1}^{n+2} = d_k + \frac{(-1)^n}{d_k^{n+1}} (\omega - d_k)^{n+1} [n \cdot \omega + d_k] \\ (\text{PPM}^L)_{n,1}^{n+2} = d_k + \frac{1}{(1-d_k)^{n+1}} (\omega - d_k)^{n+1} [n \cdot (1-\omega) + (1-d_k)] \end{cases},$$

$$\begin{cases} (\text{PPM}^L)_{n,2}^{n+3} = d_k + \frac{(-1)^n}{d_k^{n+2}} (\omega - d_k)^{n+1} \sum_{i=0}^{2} a_i^2 \omega^{2-i} d_k^i \\ (\text{PPM}^R)_{n,2}^{n+3} = d_k + \frac{1}{(1-d_k)^{n+2}} (\omega - d_k)^{n+1} \sum_{i=0}^{2} a_i^2 (1-\omega)^{2-i} (1-d_k)^i \end{cases}, \text{ and}$$

$$\begin{cases} (\text{PPM}^L)_{n,3}^{n+4} = d_k + \frac{(-1)^n}{d_k^{n+3}} (\omega - d_k)^{n+1} \sum_{i=0}^{3} a_i^3 \omega^{3-i} d_k^i \\ (\text{PPM}^R)_{n,3}^{n+4} = d_k + \frac{1}{(1-d_k)^{n+3}} (\omega - d_k)^{n+1} \sum_{i=0}^{2} a_i^3 (1-\omega)^{3-i} (1-d_k)^i \end{cases}$$

where $a_0^2 = \frac{n(n+1)}{2}$, $a_1^2 = n$ and $a_2^2 = 1$; $a_0^3 = \frac{n(n+1)(n+2)}{3!}$, $a_1^3 = \frac{n(n+1)}{2}$, $a_2^3 = n$ and $a_3^3 = 1$.

(c) $\text{IM}_{n,0;A}^{n+1}$ [8]

The mapping is a single rational polynomial which is defined as

$$\text{IM}_{n,0;A}^{n+1} = d_k + \frac{A \cdot (\omega - d_k)^{n+1}}{A \cdot (\omega - d_k)^n + \omega \cdot (1-\omega)}$$ in [0, 1] with $A$ as the free parameter and called as

IM($n$, $A$) in Ref. 8. One can check $\text{IM}_{n,0;A}^{n+1}$ satisfies $C_{n,0}$ condition. It was found by Feng et al.

[8] that $g_M = \frac{\omega \left[ \omega^2 - 3d_k \omega + (d_k + 1)d_k \right]}{(1 - 2d_k)\omega + d_k^2}$ by Henrick et a. [3] can be re-formulated in the form

of $\text{IM}_{2,0;A}^{3}$ with $n=2$ and $A=1$. $\text{IM}_{2,0;A}^{3}$ is suitable for WENO5 when critical points occur with the order 1, and the suggested choice for $A$ is $A=0.1$. As mentioned in the introduction, such choice extremely amplifies the weights at the endpoints, which was opposed in Ref. 7. However, $\text{IM}_{n,0;A}^{n+1}$ indicated a flexibility to adjust Flatness-II by choosing $A$ other than increasing the order $n$, which favors the save of computation.

(d) $\text{RM}_{n,m,k}^{n+1}$ [10]

In Ref. 10, Wang et al. proposed a specific, single rational polynomial mapping in [0, 1] as $d_k + \frac{1}{\sum_{i=0}^{m+1} a_i \omega^i} (\omega - d_k)^{n+1}$. By appealing to $C_{n,m,k}$ with $g'(0,1) = 0$, the coefficients can be obtained providing the solution exists. As the result, the meaningful formula is available only if: (i) either $m$ or $k$ should be zero, and the choice of $k=0$ is thought to render better performance. Further, (ii) $m$ and $n$ should be even. The suggested mapping is: $\text{RM}_{6,2,0}^{7} = d_k + \frac{1}{\sum_{i=0}^{3} a_i \omega^i} (\omega - d_k)^7$ with $n=6$, $m=2$ and $k=0$, where $a_0 = d_k^6$, $a_1 = -7d_k^5$, $a_2 = 21d_k^4$, $a_3 = (1-d_k)^6 - \sum_{i=0}^{2} a_i$. In the reference the mapping was called as RM(260), and $\varepsilon$ in Eq. (4) was told to take $10^{-99}$ [10]. As shown in Refs. 10-11, $\text{RM}_{6,2,0}^{7}$ suits for WENO-JS with orders up to 9, and as warned by the authors, $\text{RM}_{n,m,0}^{n+1}$ might be singular if $n$ or $m$ is odd. One can check $\text{RM}_{6,1,0}^{7}$ is ill-defined when $d_k=6/10$.

(e) $\text{AIM}_{n,m;s}^{n+1}$ [11]

Inspired by $\text{IM}_{n,0;A}^{n+1}$, Vevek et al. [11] modified the exponent of $\left[ \omega_k (1-\omega_k) \right]$ from 1 to $m+1$

and obtained: $\text{AIM}_{n,m,s}^{n+1}(\omega) = d_k + \dfrac{(\omega - d_k)^{n+1}}{(\omega - d_k)^n + s[\omega(1-\omega)]^{m+1}}$. They pointed out that $\text{AIM}_{n,m,s}^{n+1}$ would satisfy $C_{n,1}$, however one can check $\text{AIM}_{n,m,s}^{n+1}$ would actually fulfill $C_{n,m}$. As shown in Ref. 11, Flatness-II and *PEC* are competitive to each other, and therefore the ideal mapping is hard to be achievable in the case of fixed *s*. Concerning this, they further proposed an adaptive *s* as: $s = cd^{-1}\lambda$, where $\lambda = \dfrac{\min(IS_j)}{\max(IS_j) + \varepsilon_m}$ and $\varepsilon_m = \Delta x^7$. The suggested mapping for WENO7 [11] are $\text{AIM}_{4,2,1E4}^{5}$, which manifests not only good stability in long time computation but also high resolution. In Euler equations where the stencils are classified into left- and right-biased groups according to flux splitting, $\lambda$ is further modified as $\lambda = \min(\lambda_L, \lambda_R)$ with subscripts "*L, R*" denoting aforementioned group, and the scheme is referred as $\text{AIM}_{4,2,1E4}^{5}$-M in this study. However, as shown later, it is found that $\text{AIM}_{4,2,1E4}^{5}$ does not work in the case of WENO5, and even in the case of WENO7 insufficient robustness is observed.

(f) Our recent practices on piecewise rational mapping

Rather than starting with a specific form of mapping as Refs. 8 and 10-11, we try to derive piecewise rational mapping from a general formulation to satisfy $C_{n,m}$. The first step is to define a piecewise mapping at [$d_k$, 1] as $R_{n,m}^{R,n+1}(\omega) = P_{n+1}/P_{m+1}$, where $P_{n+1} = \sum\limits_{i=0}^{n+1} a_i \omega^i$ and $P_{m+1} = \sum\limits_{i=0}^{m+1} b_i \omega^i$. Without losing generality, we assume $n \geq m \geq 0$ and $a_{n+1} = 1$. It is trivial that the number of conditions in $C_{n,m}$ is $n+m+2$, while $R_{n,m}^{R,n+1}$ has $n+m+1$ unknown coefficients. Hence, there will be a free one if coefficients are solved by satisfying $C_{n,m}$. For convenience, we can take $b_{m+1}$ as the free one and drop its subscript for brevity. One can verify that meaningful solutions are achievable when $n \geq 1$. For example, that of $R_{1,0}^{R,2}$ and $R_{2,0}^{R,3}$ are $R_{1,0}^{R,2}(\omega) = \dfrac{\omega[\omega + (b-2)d_k]}{b \cdot \omega - d_k}$ and $R_{2,0}^{R,3}(\omega) = \dfrac{\omega^3 - 3d_k\omega^2 + (bd_k + 3d_k^2)\omega - bd_k - 2d_k^2 + d_k}{b\omega + (d_k - 1)^2 - b}$, which satisfy $C_{1,0}$ and $C_{2,0}$ respectively. Furthermore, it is found that they can be further rearranged as

$$R_{1,0}^{R,2}(\omega) = d_k + \frac{(\omega - d_k)^2}{-b \cdot (1-\omega) + (1-d_k)} \quad \text{and} \quad R_{2,0}^{R,3}(\omega) = d_k + \frac{(\omega - d_k)^3}{-b \cdot (1-\omega) + (1-d_k)^2}.$$ One can

see that $g_M$ can be reproduced from the latter by choosing $b = (1 - 2d_k)$. By means of similar operations, a series of $R_{n,m}^{R,n+1}(\omega)$ in $[d_k, 1]$ which satisfies $C_{n,m}$ can be derived as:

$$R_{n,m}^{R,n+1} = d_k + \frac{(\omega - d_k)^{n+1}}{c_{n,m,1}^{R}(\omega - d_k)^{n_1} + c_{n,m,2}^{R}(1-\omega)^{m+1} + c_{n,m,3}^{R}(1-d_k)^{n} + c_{n,m,4}^{R}(1-\omega)^{m}} \quad (11)$$

where $n_1$, coefficients $c_{n,m,i}^{R}$ except $c_{n,m,4}^{R}$ are tabulated in Table 1 and 2 for $m, n \leq 4$. Almost all $c_{n,m,4}^{R}$ equals zero except $c_{4,2,4}^{R} = 4(d_k - 1)^2$.

Table 1 $n_1$ of $(n, m)$ in mapping $R_{n,m}^{R,n+1}$

|     | m=0 | m=1 | m=2 | m=3 | m=4 |
|-----|-----|-----|-----|-----|-----|
| n=1 | 0   | 1   | -   | -   | -   |
| n=2 | 0   | 1   | 2   | -   | -   |
| n=3 | 0   | 1   | 3   | 3   |     |
| n=4 | 0   | 1   | 2   | 4   | 4   |

Table 2 Coefficients $c_{n,m,i}^{R}$ except $c_{n,m,4}^{R}$ of mapping $R_{n,m}^{R,n+1}$

|     |     | m=0 | m=1 | m=2 | m=3 | m=4 |
|-----|-----|-----|-----|-----|-----|-----|
|     | i=1 | 0   | 1   | -   | -   | -   |
| n=1 | i=2 | -b  | $\frac{b}{1-b \cdot d_k}$ | - | - | - |
|     | i=3 | 1   | 0   | -   | -   | -   |
|     | i=1 | 0   | $2(1-d_k)$ | 1 | - | - |
| n=2 | i=2 | -b  | b   | $\frac{-b}{1-b \cdot d_k}$ | - | - |
|     | i=3 | 1   | -1  | 0   | -   | -   |
|     | i=1 | 0   | $3(1-d_k)^2$ | 1 | 1 | - |
| n=3 | i=2 | -b  | b   | 1-b | $\frac{b}{1-b \cdot d_k}$ | - |
|     | i=3 | 1   | -2  | 0   | 0   | -   |
|     | i=1 | 0   | $4(1-d_k)^3$ | $2(d_k-1)^2$ | 1 | 1 |
| n=4 | i=2 | -b  | b   | -b  | -(1-b) | $\frac{-b}{1-b \cdot d_k}$ |
|     | i=3 | 1   | -3  | -1  | 0   | 0   |

Theoretically, $R_{n,m}^{L,n+1}$ in [0, $d_k$] satisfying $C_{n,m}$ can be derived similarly. However, on observing its symmetry with $R_{n,m}^{R,n+1}$ about $\omega = d_k$, $R_{n,m}^{L,n+1}$ can be obtained through the following transformation:

(i) $\mathbf{F}^{(1)}(\omega) = R_{n,m}^{R,n+1}(\omega + d_k) - d_k$, where $\omega \in [0, 1-d_k]$ and $\mathbf{F}^{(1)} \in [0, 1-d_k]$,

(ii) $\mathbf{F}^{(2)}(\omega) = \frac{d_k}{1-d_k} \mathbf{F}^{(1)}(\frac{1-d_k}{d_k}\omega)$, where $\omega \in [0, d_k]$ and $\mathbf{F}^{(2)} \in [0, d_k]$,

(iii) $\mathbf{F}^{(3)}(\omega) = -\mathbf{F}^{(2)}(-\omega)$, where $\omega \in [-d_k, 0]$ and $\mathbf{F}^{(3)} \in [-d_k, 0]$,

(iv) Finally, $R_{n,m}^{L,n+1} = \mathbf{F}^{(3)}(\omega - d_k) + d_k$ where $\omega \in [0, d_k]$ and $R_{n,m}^{L,n+1} \in [0, d_k]$.

The above procedure can be further assembled as: $g^L = \frac{d_k}{1-d_k}\left[1 - g^R(1 - \frac{1-d_k}{d_k}\omega)\right]$ where $g^R$ and $g^L$ corresponds to $R_{n,m}^{R,n+1}$ and $R_{n,m}^{L,n+1}$ here. One can verify that the acquired mapping would satisfy $C_{n,m}$ providing $R_{n,m}^{L,n+1}$ satisfies the same condition. Finally, similar $R_{n,m}^{L,n+1}$ can be derived as

$$R_{n,m}^{L,n+1} = d_k + \frac{(\omega - d_k)^{n+1}}{c_{n,m,1}^L(\omega - d_k)^{n_1} + c_{n,m,2}^L\omega^{m+1} + c_{n,m,3}^L d_k^{\ n} + c_{n,m,4}^L\omega^m} \quad (12)$$

where coefficients $c_i^{L,n}$ tabulated in Table 13 in Appendix II for $m, n \leq 4$ for completeness.

Regarding $R_{n,m}^{R/L,n+1}$, the following remarks are given:

(i) Although there exists one free parameter $b$, it cannot further make either $C_{n+1,m}$ or $C_{n,m+1}$ achievable. One can verify the requirement of $C_{n+1,m}$ or $C_{n,m+1}$ will render $R_{n,m}^{R,n+1}$ insolvable, singular or reduced to $R_{n,m}^{R,n+1} = d_k$ or $\omega$.

(ii) In order to avoid zero point(s) in the denominator of Eq. (11) in [$d_k$, 1], the range of $b$ should be confined. Additionally, on making the mapping lay below the identity mapping, the relation $(-1)^m \frac{\partial^{m+1} R_{n,m}^{R,n+1}}{\partial x^{m+1}}(1) > 0$ should be required, through which additional confinement of $b$ can be acquired also. Considering the two confinements, the range of $b$ for $R_{n,m}^{R,n+1}$ can be obtained and tabulated in Table 3. One can testify $R_{n,m}^{R,n+1}$ is of singularity and lay below the identity mapping when $b$ is within the range. As expected, the same range of $b$ exists for $R_{n,m}^{L,n+1}$ also.

Table 3 Valid range of $b$ for $R_{n,m}^{R,n+1}$

|     | m=0 | m=1 | m=2 | m=3 | m=4 |
| --- | --- | --- | --- | --- | --- |
| n=1 | $b<1$ | $0<b<1/d_k$ | - | - | - |
| n=2 | $b<1-d_k$ | $b>1$ | $b<0 \cup b>\frac{1}{d_k}$ | - | - |
| n=3 | $b<(1-d_k)^2$ | $b\geq 3(1-d_k)$ | $b<1$ | $0<b<1/d_k$ | - |
| n=4 | $b<(1-d_k)^3$ | $b\geq 6(1-d_k)^2$ | $b<3(1-d_k)$ | $b>1$ | $b<0 \cup b>\frac{1}{d_k}$ |

Combining Table 3 and 2, one can get another representation of $c_{n,m,i}^R$ as shown in Table 4, where $c_{n,m,2}^R$ acts as the free parameter.

Table 4 Another recipe of $c_{n,m,i}^R$ except $c_{n,m,4}^R$ of mapping $R_{n,m}^{R,n+1}$

|     |     | m=0 | m=1 | m=2 | m=3 | m=4 |
| --- | --- | --- | --- | --- | --- | --- |
| n=1 | i=1 | 0 | 1 | - | - | - |
|     | i=2 | $c_{n,m,2}^R > -1$ | $c_{n,m,2}^R > 0$ | - | - | - |
|     | i=3 | 1 | 0 | - | - | - |
| n=2 | i=1 | 0 | $2(1-d_k)$ | 1 | - | - |
|     | i=2 | $c_{n,m,2}^R > -(1-d_k)$ | $c_{n,m,2}^R > 1$ | $c_{n,m,2}^R > 0$ | - | - |
|     | i=3 | 1 | -1 | 0 | - | - |
| n=3 | i=1 | 0 | $3(1-d_k)^2$ | 1 | 1 | - |
|     | i=2 | $c_{n,m,2}^R > -(1-d_k)^2$ | $c_{n,m,2}^R \geq 3(1-d_k)$ | $c_{n,m,2}^R > 0$ | $c_{n,m,2}^R > 0$ | - |
|     | i=3 | 1 | -2 | 0 | 0 | - |
| n=4 | i=1 | 0 | $4(1-d_k)^3$ | $2(1-d_k)^2$ | 1 | 1 |
|     | i=2 | $c_{n,m,2}^R > -(1-d_k)^3$ | $c_{n,m,2}^R \geq 3(1-d_k)^2$ | $c_{n,m,2}^R > -3(1-d_k)$ | $c_{n,m,2}^R > 0$ | $c_{n,m,2}^R > 0$ |
|     | i=3 | 1 | -3 | -1 | 0 | 0 |

It is noted from the table the valid range of $c_{n,n,2}^R$ is: $c_{n,n,2}^R > 0$, and it is trivial that the denominator in $R_{n,n}^{R,n+1}$, namely $(\omega-d_k)^n + c(1-\omega)^{m+1}$, is free of zero point at $\omega \in [d_k, 1]$ providing $c>0$. Motivated by this observation, the following theorem is proposed:

**Theorem 2**: Consider a mapping as $R_{n,m}^{n+1} = d_k + \dfrac{(\omega-d_k)^{n+1}}{(\omega-d_k)^n + c(1-\omega)^{m+1}}$ in $[d_k, 1]$ with $0<d_k<1$ and $n,m \geq 1$. $c>0$ is the sufficient and necessary condition for $R_{n,m}^{n+1}$ to be void of singularity.

For brevity, the proof of the theorem is given in Appendix III.

One can verify similar conclusion be established, i.e. the mapping $d_k + \dfrac{(\omega - d_k)^{n+1}}{(\omega - d_k)^n + c\omega^{m+1}}$ at $\omega \in [0, d_k]$ with $0 < d_k < 1$ and $n, m \geq 1$ will be free of singularity providing $(-1)^n c > 0$.

(iii) Consider $R_{n,m}^{R,n+1}$ by Eq. (11). $\partial^{n+1} R_{n,m}^{R,n+1}(d_k)/\partial \omega^{n+1}$ and $\partial^{m+1} R_{n,m}^{R,n+1}(1)/\partial \omega^{m+1}$ can be derived and tabulated in Table 5.

Table 5 $R_{n,m}^{R,n+1(n+1)}(d_k)$ and $R_{n,m}^{R,n+1(m+1)}(1)$ where $d_{1k}$ denotes $(1 - d_k)$

|  | m=0 | m=1 | m=2 | m=3 | m=4 |
|---|---|---|---|---|---|
| $R_{1,m}^{R,2(2)}(d_k)$ | $\dfrac{2}{(c_{n,m,2}^R + 1)d_{1k}}$ | $\dfrac{2}{c_{n,m,2}^R d_{1k}^2}$ | - | - | - |
| $R_{1,m}^{R,2(m+1)}(1)$ | $c_{n,m,2}^R + 2$ | $-2c_{n,m,2}^R$ | - | - | - |
| $R_{2,m}^{R,3(3)}(d_k)$ | $\dfrac{6}{c_{n,m,2}^R d_{1k} + d_{1k}^2}$ | $\dfrac{6}{c_{n,m,2}^R d_{1k}^2 - d_{1k}^2}$ | $\dfrac{6}{c_{n,m,2}^R d_{1k}^3}$ | - | - |
| $R_{2,m}^{R,3(m+1)}(1)$ | $\dfrac{c_{n,m,2}^R + 3d_{1k}}{d_{1k}}$ | $\dfrac{-2c_{n,m,2}^R + 2}{d_{1k}}$ | $\dfrac{6c_{n,m,2}^R}{d_{1k}}$ | - | - |
| $R_{3,m}^{R,4(4)}(d_k)$ | $\dfrac{24}{c_{n,m,2}^R d_{1k} + d_{1k}^3}$ | $\dfrac{24}{c_{n,m,2}^R d_{1k}^2 - 2d_{1k}^3}$ | $\dfrac{24}{c_{n,m,2}^R d_{1k}^3}$ | $\dfrac{24}{c_{n,m,2}^R d_{1k}^4}$ | - |
| $R_{3,m}^{R,4(m+1)}(1)$ | $\dfrac{c_{n,m,2}^R + 4d_{1k}^2}{d_{1k}^2}$ | $\dfrac{-2c_{n,m,2}^R + 6d_{1k}}{d_{1k}^2}$ | $\dfrac{6c_{n,m,2}^R}{d_{1k}^2}$ | $-\dfrac{24c_{n,m,2}^R}{d_{1k}^2}$ | - |
| $R_{4,m}^{R,5(5)}(d_k)$ | $\dfrac{120}{c_{n,m,2}^R d_{1k} + d_{1k}^4}$ | $\dfrac{120}{c_{n,m,2}^R d_{1k}^2 - 3d_{1k}^4}$ | $\dfrac{120}{c_{n,m,2}^R d_{1k}^3 + 3d_{1k}^4}$ | $\dfrac{120}{c_{n,m,2}^R d_{1k}^4}$ | $\dfrac{120}{c_{n,m,2}^R d_{1k}^5}$ |
| $R_{4,m}^{R,5(m+1)}(1)$ | $\dfrac{c_{n,m,2}^R + 5d_{1k}^3}{d_{1k}^3}$ | $\dfrac{-2c_{n,m,2}^R + 12d_{1k}^2}{d_{1k}^3}$ | $\dfrac{6c_{n,m,2}^R + 12d_{1k}}{d_{1k}^3}$ | $-\dfrac{24c_{n,m,2}^R}{d_{1k}^3}$ | $\dfrac{120c_{n,m,2}^R}{d_{1k}^3}$ |

From the table, one can see that $R_{n,m}^{R,n+1(n+1)}(d_k)$ will decrease with the increase of $c_{n,m,2}^R$, which indicates the mapping will become more flat as away from $\omega = d_k$, whereas $(-1)^m R_{n,m}^{R,n+1(m+1)}(1)$ will increase which indicates the mapping will diverge more from the identity mapping as away from $\omega = 1$; and the converse is also true. Hence, Flatness-II and $PEC$ of $R_{n,m}^{R,n+1}$ would be competitive with each other, or they can hardly reach to a desired state at the same time.

(iv) Recalling previous formulation of $\text{AIM}_{n,m;s}^{n+1}$, one can find that $R_{n,n}^{R(L),n+1}$ quite

resembles $\text{AIM}_{n,n;s}^{n+1}(\omega)$ if the overall term $s[\omega(1-\omega)]^{m+1}$ in the denominator can be somehow re-interpreted as $(1-\omega)^{m+1}$ in $[d_k, 1]$ and $\omega^{m+1}$ in $[0, d_k]$. As just analyzed, the fixed value of $s$ would constrain the performance of mapping, therefore we propose a new formulation other than the adaptive manner in $\text{AIM}_{n,n;s}^{n+1}(\omega)$ but with sufficient controllability in Section 4. Prior to further investigations, several uncertainties of mappings will be investigated.

### 3 Investigation on uncertainties of mapping methods

As referred in introduction, several uncertainties regarding mapping exist, and the investigation of which may help to further develop mapping methods.

(1) Behavior of endpoint convergence

In Refs. 7 and 10, *PEC* of mapping was suggested to ENO by appealing to $g^{(i)}(0,1) = 0$ with $i \leq m$ for given $m$. In the reference, the long-time simulations of scalar advection of square wave and combination-waves were studied, where WENO5-JS was tested for reference which indicates large errors. Such consequences might lead to an impression that *PEC* toward WENO would possibly cause large errors. Because alternative *PEC* as WENO was suggested by Ref. 9, it is necessary to clarify if such practice could yield results similar to that by *PEC* as ENO. For this purpose, the following numerical study is planned: first, WENO5-JS is still employed as the benchmark scheme as in Refs. 7 and 10; next $PM^8_{6,1}$ is chosen as the representative with *PEC* as ENO, and for comparison the equivalent $PPM^8_{6,1}$ is used with *PEC* as WENO; thirdly, the scalar advection of combination-waves at grid points $N=800$ are tested with computation time $t=2000$, which was regarded as the typical case to check stability in long time computation. Corresponding results are shown in Fig. 1, where $g_M$ is tested also. From the figure, $PPM^8_{6,1}$ yields a result which is very similar to that by $PM^8_{6,1}$ [8], where nearly all distribution is smooth except little oscillations near $x=0.6$. Hence, *PEC* as WENO indicates a comparable performance and therefore the same applicability as *PEC* as ENO on developing mapping function.

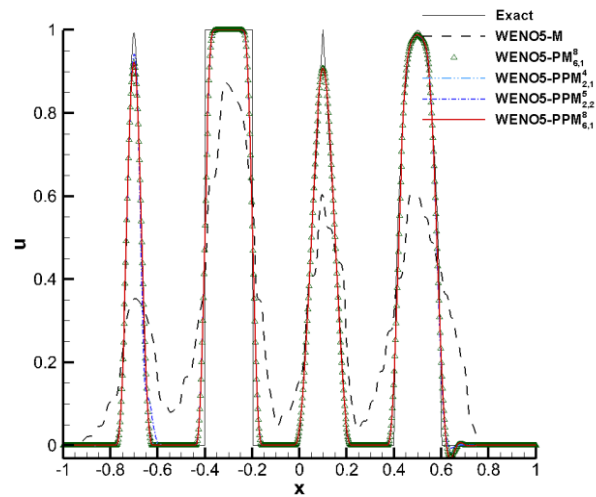

(a) Overall view

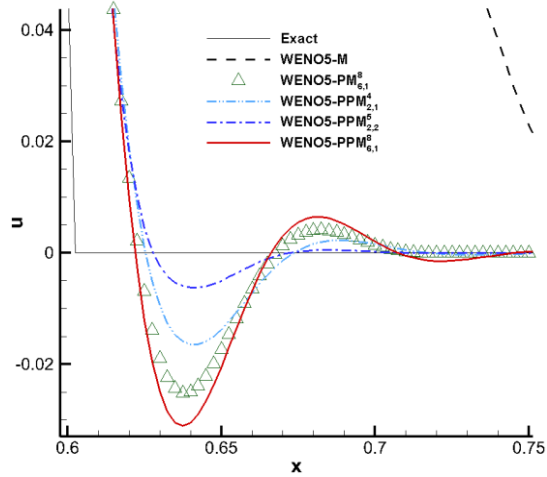

(b) Zoomed view around $x=0.6$

Fig. 1. Results of scalar advection of combination-waves by WENO5-PM$^8_{6,1}$, -PPM$^8_{6,1}$, -PPM$^5_{5,2}$, -PPM$^4_{2,1}$, and WENO5-M with $N=800$ and $t=2000$

Besides the use of PM$^8_{6,1}$ and PPM$^8_{6,1}$ for WENO5-JS, Section 4 and Table 12 in Appendix I will show that the minimum "$n$" in $C_{n,m}$ for optimal order recovery therein would be 2 other than 6, hence corresponding PPMs, e.g. PPM$^4_{2,1}$ and PPM$^5_{2,2}$ can be applied as well. For reference, they are also tested and yield results shown in zoomed view of Fig. 1(b). The performances indicate that oscillations by both mappings decrease, especially the result of PPM$^5_{2,2}$ only shows tiny undershoot. The consequences confirm the capability and well attributes of *PEC* as WENO.

If further comparing the results of above PPM$^{6+m+1}_{6,m}$ and PPM$^{2+m+1}_{2,m}$, one can find that the former is relatively more oscillatory than that of the latter. The possible cause is that with the increase of "$n$" in $C_{n,m}$, both Flatness-I and –II increase, which is liable to numerical instability. Hence in order to recover optimal order at critical point, the employment of larger "$n$" in $C_{n,m}$ than necessary is not suggested unless really needed.

(2) Effect of piecewise implementation on numerical stability

In Ref. 10, Wang et al. suggested the property of smoothness for a good mapping and expected the function to have infinite derivatives in [0, 1]. It is clear that the single continuous function would fulfill the expectation. As an illustration, the computation of scalar advection of combination-waves was made by WENO9-PM$^8_{6,1}$ and oscillatory result was obtained. They concluded the oscillations owed to the insufficient smoothness of PM$^8_{6,1}$, or its piecewise implementation. Because various possibilities that affect numerical stability might exist, the above conclusion deserves further investigation. In this regard, the similar commensurate PPM$^8_{6,1}$ is employed with WENO9 in the same problem, where 400 grid points are used and the computation advances to $t=100$ as in Ref. 10. From the figure, PPM$^8_{6,1}$ indicates an oscillatory distribution almost the same as that of PM$^8_{6,1}$, which seems to favor the conclusion of Ref. 10. However, disparate results emerge when more mappings are included for testing as follows.

Two group of piecewise mappings are chosen: the first group is PPM$^6_{4,1}$, PPM$^7_{4,2}$ and PPM$^8_{4,3}$ where $n=4$ for $C_{n,m}$ plus PM$^8_{6,1}$, and the second group is PPM$^8_{6,1}$, PPM$^9_{6,2}$, PPM$^{10}_{6,3}$ and PM$^8_{6,1}$ with $n=6$. Corresponding results are shown in Figs. 2(b)-(d). In the first group, on the one hand the results of PPM$^6_{4,m}$ show less oscillations than that of PM$^8_{6,1}$ and PPM$^8_{6,1}$, on the other hand oscillations decrease with the increase of *m* and are almost disappeared when *m*=3. In the second group, all mappings with the same *m* from 1-3 yield oscillatory results. Hence the oscillations by WENO9-JS with piecewise mapping cannot be attributed to insufficient smoothness by piecewise function as claimed by Ref. 10, but come from the inadequate numerical stability caused by too flat profile around $\omega = d_k$ and which does not make the mapping tend to identity mapping fast enough. Based on above observations, the larger value of *m* in $C_{n,m}$ to define mapping, e.g. *m*=*n*, is suggested in subsequent studies.

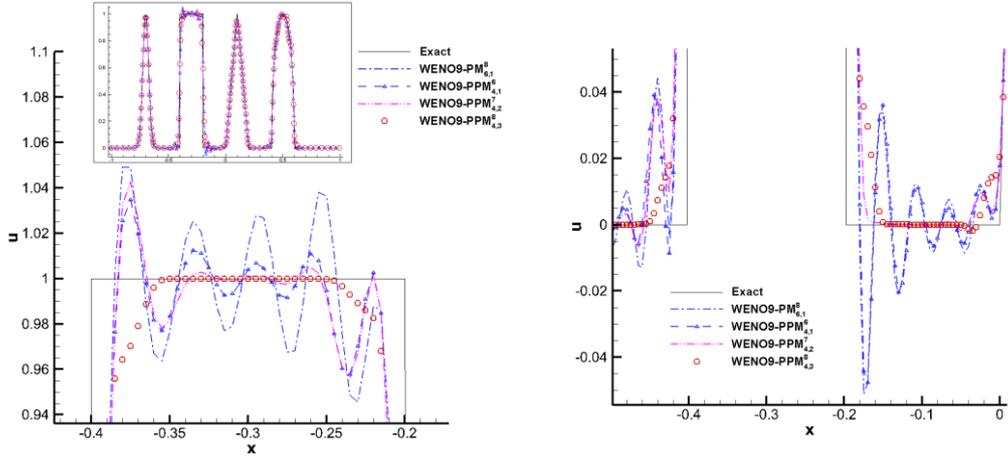

(a) Zoomed view of results of first group around the peak of the square wave

(b) Zoomed view of results of first group around the foot of the square wave

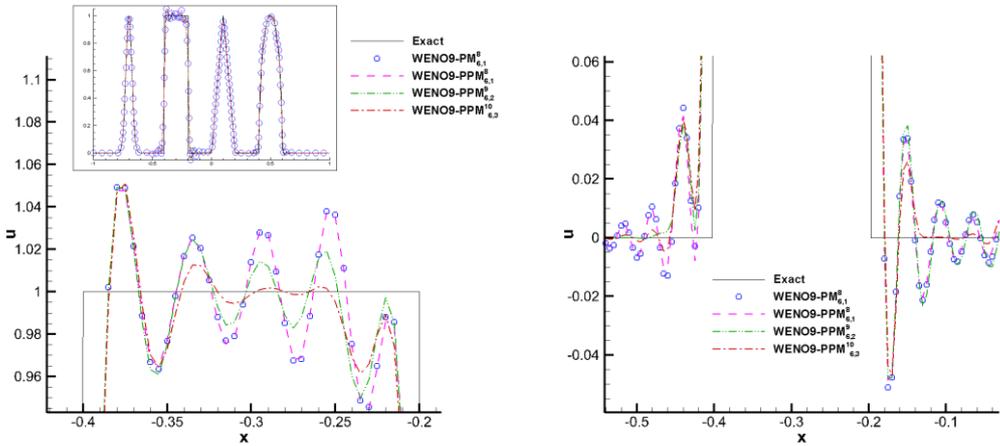

(c) Zoomed view of results of second group around the peak of the square wave

(d) Zoomed view of results of second group around the foot of the square wave

Fig. 2. Result comparison between the first group {WENO9-PPM$^6_{4,1}$, -PPM$^7_{4,2}$, -PPM$^8_{4,3}$, -PM$^8_{6,1}$} and the second group {WENO9-PM$^8_{6,1}$, WENO9-PM$^8_{6,1}$, -PPM$^8_{6,1}$, -PPM$^{10}_{6,3}$} on scalar advection of combination-waves with *N*=400 and *t*=100

To further validate the above statement, the case of mapped WENO-JS is investigated under *N*=200 and *t*=1000 where PPM$^5_{3,1}$, PM$^8_{6,1}$ and PPM$^8_{6,1}$ are employed. is provide in Fig. 3. One can see that derivatives of PPM$^5_{3,1}$ at $\omega_k = d_k$ from both sides only equal to each other up to the third-order, which is much smaller than the sixth-order of PM$^8_{6,1}$. In the view of Ref. 10,

PPM$^5_{3,1}$ appears unsmooth in the sense of continuity, which would entail prominent oscillations very likely. However, Fig. 3 tells that PPM$^5_{3,1}$ yields a quite smooth result as that of PM$^8_{6,1}$, which indicates the discontinuous fourth-order derivatives would not cause numerical instability. One may observe the error in the corner of $x$=-0.8, and according to our experience and understanding, such phenomenon owes to the dissipation of the mapping rather than the occurrence of oscillation.

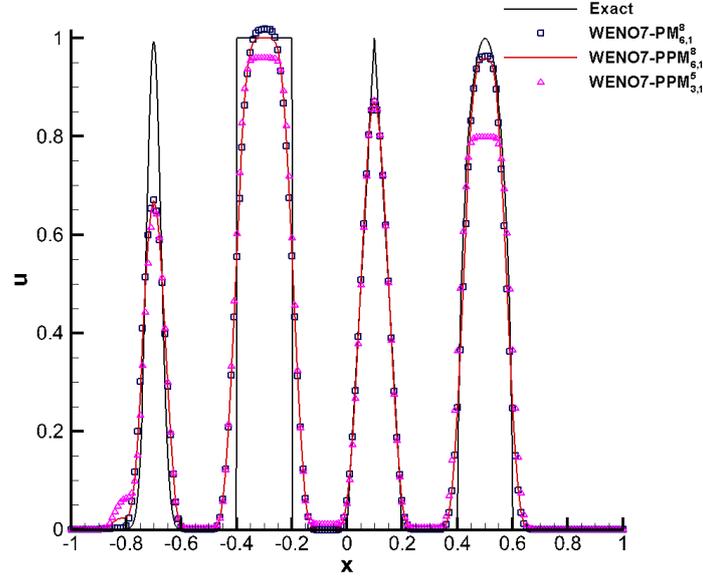

Fig. 3. Results of WENO7-PPM$^5_{3,1}$, PM$^8_{6,1}$ and PPM$^8_{6,1}$ on scalar advection of combination-waves with grid points $N$=200 and computation time $t$=1000

Hence, the piecewise implementation of mapping doesn't necessarily entail numerical stability.

(3) Feasibility to recover optimal order of WENO3-JS at the first-order critical point

Refs. 3, 8 and 10 showed the incapability of mappings to make WENO3-JS preserve the optimal order in the occurrence of first-order critical points. The mechanism is regarded as: when the critical point with the order $n_{cp}$ occurs, $\omega_k - d_k = O(\Delta x^{r-1-n_{cp}})$ where $r$ is the order of substencil and $n_{cp}$ indicates $f' = ... = f^{(n_{cp})} = 0$ and $f^{(n_{cp}+1)} \neq 0$. In the case of WENO3-JS with $r$=2 and when $n_{cp}$=1 occurs, $\omega_k - d_k = O(1)$ and therefore the least condition required by mapping is violated. Considering the situation absent of critical points, the smoothness indicators have accuracy relations as $IS_k^2 = (f'_j \Delta x)^2 + (-1)^{k+1} f'_j f''_j \Delta x^3 + \left(\frac{1}{3} f'_j f'''_j + \frac{1}{4} f''^2_j \right) \Delta x^4 + O(\Delta x^5)$. However, when the first-order critical point occurs at $x_j$, one can find that $IS_k^2 = \frac{1}{4}(f''_j \Delta x^2)^2 \left[1 + (-1)^{k+1} \frac{2}{3} \frac{f'''_j}{f''_j} \Delta x + O(\Delta x^2)\right]$. So it seems that WENO3-JS can still preserve the optimal order though mapping, and in fact such preservation does not happen in tests such as 1D scalar advection with the initial distribution $sin(x)$.

A heuristic analysis in this regard is proposed as: when the first-order critical point occurs at $x_j$, the position will inevitably shift to location such as $x_c \in [x_j, x_{j+1}]$ next. For convenience, suppose $x_j - x_c = \frac{-1}{m} \Delta x$ where $1 < m < \infty$. Hence $f'_j = f''_{x_c}(x_j - x_c) + O(\Delta x^2)$, and it can be derived

that

$$IS_0^2 = f''^2_{x_c}(x_j - x_c)^2 \Delta x^2 - (x_j - x_c) f''^2_{x_c} \Delta x^3 + \frac{1}{4} f''^2_{x_c} \Delta x^4 + O(\Delta x^5) = (\tfrac{1+m}{m^2} + \tfrac{1}{4}) f''^2_{x_c} \Delta x^4 + O(\Delta x^5)$$

and

$$IS_1^2 = f''^2_{x_c}(x_j - x_c)^2 \Delta x^2 + (x_j - x_c) f''^2_{x_c} \Delta x^3 + \frac{1}{4} f''^2_{x_c} \Delta x^4 + O(\Delta x^5) = (\tfrac{1-m}{m^2} + \tfrac{1}{4}) f''^2_{x_c} \Delta x^4 + O(\Delta x^5)$$

. Hence $IS_k = D(1 + O(\Delta x))$ and consequent $\omega_k - d_k = O(1)$ are not established, thereby the optimal order would not be preserved.

One of remedies for the problem is to upgrade smoothness indicators of WENO3-JS to that of WENO5-JS as suggested by Ref. 14. Specifically, $IS_0^2$ and $IS_1^2$ are substituted by $IS_0^3$ and $IS_2^3$ respectively. In the case of the first-order critical point occurs at $x_j$, it is well-known that $IS_k^3 = \frac{13}{12}(f''_j \Delta x^2)^2 (1 + O(\Delta x))$ [3] and therefore the recovery of optimal order is achievable. When the critical point occurs at aforementioned $x_c$, one can check that $IS_k^3 = (\frac{1}{m^2} + \frac{13}{12})(f''_j \Delta x^2)^2 (1 + O(\Delta x))$, and the least requirement of mapping is still satisfied.

In short, the order of WENO3-JS by mapping at the occurrence of first-order critical point is available through upgrading smoothness indicator.

**4 New piecewise rational mapping functions**

As discussed in Section 2.2, for existing mappings in fixed other than dynamic form, the achievement of flatness and *PEC* at the same time is competitive and thereby unachievable. This consequence indicate the increase of flatness would suppress the endpoint convergence, which potentially impair the numerical stability. In the mean while recent studies show more and more concerns on the stability in long time computation [7,8,10-11], and it has been found that although various mappings usually have resolutions improved in short time, they might perform oscillatory in the long run. According to pervious discussion, the mechanism is suspected to be the mismatch between characteristics such as flatness and endpoint convergence. In this section, we propose a new method with sufficient regulation which would make the mapping satisfy prescribed $C_{n,m}$ on the one hand, on the other hand could achieve desirable integration of the flatness and *PEC*.

Based on our recent experiences in point (f) of Section 2.2 while observing the practice of $\text{AIM}_{n,m,s}^{n+1}$ [11], a new piecewise rational mapping PRM in $[d_k, 1]$ is proposed as:

$$\text{PRM}_{n,m;n_1,m_1;c_1,c_2}^{R,n+1} = d_k + \frac{(\omega - d_k)^{n+1}}{(\omega - d_k)^n + c_2(\omega - d_k)^{n_1}(1-\omega)^{m_1} + c_1(1-\omega)^{m+1}} \tag{13}$$

where all exponents are greater than or equal to 1, and where $m \leq n$ is suggested for practicability. As shown later, Eq. (13) satisfies $C_{n,m}$ providing $m_1 \geq m+1$. By aforementioned transformation, i.e. $g^L = \frac{d_k}{1-d_k}\left[1 - g^R(1 - \frac{1-d_k}{d_k}\omega)\right]$, corresponding PRM in $[1, d_k]$ can be obtained as:

$$\text{PRM}^{L,n+1}_{n,m;n_1,m_1;c_1,c_2} = d_k + \frac{(\omega-d_k)^{n+1}}{(\omega-d_k)^n + (-1)^{n_1+n}\left(\frac{1-d_k}{d_k}\right)^{n_1+m_1-n} c_2(\omega-d_k)^{n_1}\omega^{m_1} + (-1)^n\left(\frac{1-d_k}{d_k}\right)^{m-n+1} c_1\omega^{m+1}} \quad (14)$$

which satisfies $C_{n,m}$ as well if $m_1 \geq m+1$. Especially, if $n_1=1$, PRMs become:

$$\begin{cases} \text{PRM}^{L,n+1}_{n,m;m_1;c_1,c_2} = d_k + \dfrac{(\omega-d_k)^{n+1}}{(\omega-d_k)^n + (-1)^{1+n}\left(\frac{1-d_k}{d_k}\right)^{1+m_1-n} c_2(\omega-d_k)\omega^{m_1} + (-1)^n\left(\frac{1-d_k}{d_k}\right) c_1\omega^{m+1}} \\ \text{PRM}^{R,n+1}_{n,m;m_1;c_1,c_2} = d_k + \dfrac{(\omega-d_k)^{n+1}}{(\omega-d_k)^n + c_2(\omega-d_k)(1-\omega)^{m_1} + c_1(1-\omega)^{m+1}} \end{cases}.$$

It is conceivable that $c_i$ can take different values piecewisely, and the coefficients regarding $d_k$ in $\text{PRM}^{L/R}$ can be absorbed into $c_i$ to yield the following formulation at $n_1=1$:

$$\begin{cases} \text{PRM}^{L,n+1}_{n,m;m_1;c_1,c_2} = d_k + \dfrac{(\omega-d_k)^{n+1}}{(\omega-d_k)^n + (-1)^{1+n} c_2^L(\omega-d_k)\omega^{m_1} + (-1)^n c_1^L\omega^{m+1}} \\ \text{PRM}^{R,n+1}_{n,m;m_1;c_1,c_2} = d_k + \dfrac{(\omega-d_k)^{n+1}}{(\omega-d_k)^n + c_2^R(\omega-d_k)(1-\omega)^{m_1} + c_1^R(1-\omega)^{m+1}} \end{cases}. \quad (15)$$

When $m=n$, the equation becomes

$$\begin{cases} \text{PRM}^{L,n+1}_{n,n;m_1;c_1^L,c_2^L} = d_k + \dfrac{(\omega-d_k)^{n+1}}{(\omega-d_k)^n + (-1)^{1+n} c_2^L(\omega-d_k)\omega^{m_1} + (-1)^n c_1^L\omega^{n+1}} \\ \text{PRM}^{R,n+1}_{n,n;m_1;c_1^R,c_2^R} = d_k + \dfrac{(\omega-d_k)^{n+1}}{(\omega-d_k)^n + c_2^R(\omega-d_k)(1-\omega)^{m_1} + c_1^R(1-\omega)^{n+1}} \end{cases}. \quad (16)$$

According to Theorem 1 and the characteristics of aforementioned transformation, it is conceivable that positive $c_i^{L/R}$ would make the mappings free of singularity. Furthermore, the relationship of $\text{PRM}^{L,n+1}_{n,m;m_1;c_1,c_2}$ with $C_{n,m}$ is indicated by the following theorem with the proof shown in Appendix III.

**Theorem 3:** Consider a function as $f(\omega) = d_k + \dfrac{(\omega-d_k)^{n+1}}{(\omega-d_k)^n + c_2(\omega-d_k)^{n_1}(1-\omega)^{m_1} + c_1(1-\omega)^{m+1}}$ in $[d_k, 1]$ where $n \geq 1$, and $c_1$, $c_2$ have the same sign. If $m \geq 0$ and $m_1 \geq 1$, then $f(\omega)$ satisfies $C_{n,\min(m,m_1-1)}$.

In purpose of application, it is important to explore the functions of $c_1$, $c_2$, $n_1$, $m_1$ in $\text{PRM}^{L/R}$. In this regard, parametric studies of $\text{PRM}^{R,n+1}_{n,m;n_1,m_1;c_1,c_2}$ is made first; according to similarity, that of $\text{PRM}^{L,n+1}_{n,m;n_1,m_1;c_1,c_2}$ can be understood as well. Specifically, $\text{PRM}^{R,3}_{2,2;n_1,m_1;c_1,c_2}$ and its profile at $d_k = 3/10$ are chosen for illustration. For convenience, the whole or part of super- and/or sub-scripts of PRM and $c_i^{L/R}$ are omitted sometimes for convenience.

First, the effect of $c_1$ is discussed. As just mentioned, $c_1$ should be positive in order to void singularity of mapping. A series of $c_1$ are taken as $\{1, 10, 100\}$ under $c_2=10$, $n_1=1$ and $m_1=5$ and corresponding profiles are drawn in Fig. 4(a). One can see the only decrease of $c_1$ would yield the less flatness of mapping around $\omega=d_k$ and the increased convergence to identity mapping near $\omega=1$, and vice versa. Next, the effect of $c_2$ is studied by taking values $\{1, 10^3, 10^6\}$ under the choice $c_1=1$, $n_1=1$ and $m_1=5$, and the result is shown in Fig. 4(b). The figure tells the increase of $c_2$ will extend the flatness and accelerate endpoint convergence by pushing the junction upward, and vice versa. Thirdly, the effect of $n_1$ is explored by taking $\{1,4,10\}$ under $c_1=1$, $c_2=100$ and $m_1=2$, and corresponding results are shown in Fig. 4(c). One can see that the enlargement of $n_1$ will decrease the flatness of mapping without largely shrinking endpoint convergence. Base on this observation and in save of cost, $n_1$ is suggested to take value of 1. At last, the effect of $m_1$ is investigated by taking $\{2, 3, 5\}$ under $c_1=1$, $c_2=100$ and $n_1=1$. Theoretically $m_1$ controls how fast the influence of $c_2(\omega-d_k)^{n_1}(1-\omega)^{m_1}$ will vanish near $\omega=1$. One can observe that on the one hand $m_1$ can adjust the flatness and degree of *PEC* (i.e. the smaller the flatter), on the other hand it controls the transition slope of profile from $\omega=d_k$ to $\omega=1$. The larger $m_1$ is, the more abrupt the transition of mapping will perform.

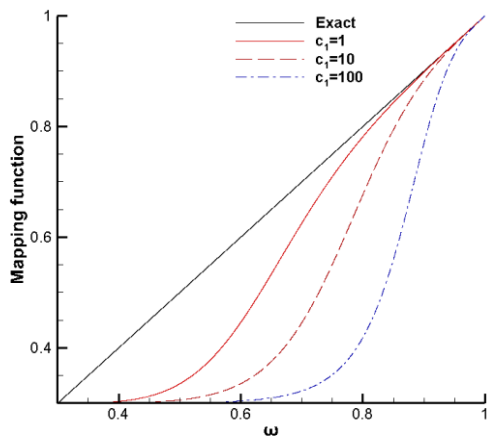
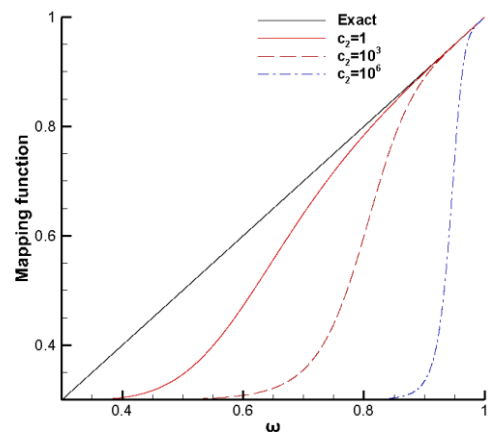

(a) Effects of $c_1$ by taking values $\{1,10,100\}$ under $c_2=10$, $n_1=1$ and $m_1=5$

(b) Effects of $c_2$ by taking values $\{1,10^3,10^6\}$ under $c_1=1$, $n_1=1$ and $m_1=5$

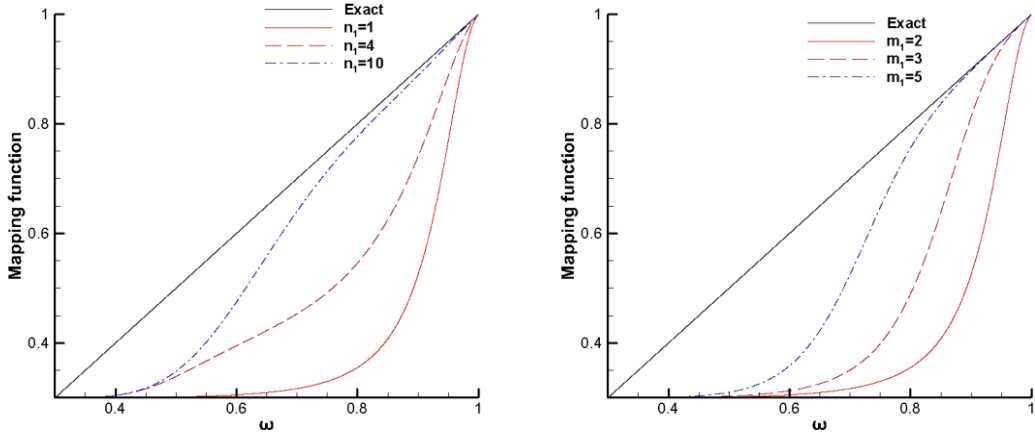

(c) Effects of $n_1$ by taking values $\{1,4,10\}$ under $c_1=1$, $c_2=100$ and $m_1=2$

(d) Effects of $m_1$ by taking values $\{2,3,5\}$ under $c_1=1$, $c_2=100$ and $n_1=1$

Fig. 4. Effects of $c_1$, $c_2$, $n_1$, $m_1$ on the profile of $\text{PRM}_{2,2}^{R,3}$ in the case of $d_k=3/10$ (solid, dash and dash-dot lines correspond to parameter values from small to big)

Base on the above discussion, extensive numerical practices are made afterwards, and the following recommendations are provided to determine parameters including $n$ and $m$:

(1) The determination of $n$

According to understandings in Refs. 3, 8 and 10, given WENO$k$-JS, the minimum required order $r_{c,g}$ of $g(\omega_k)$ such that $g^{(i)}(\omega_k)=0$ for $0 < i \le r_{c,g}$ to recover the optimal order at critical point is tabulated in Table 12 in Appendix I. In the meanwhile, $r_{c,g}$ defines the extent of Flatness-I. It is obvious that $r_{c,g}$ equals to $n$ in $C_{n,m}$. Hence for WENO3-JS, $n$ in Eqns. (13-14) can take 1; for WENO5- and WENO7-JS, $n=4$ and 5 and suggested respectively.

(2) The determination of $m$

As $m$ corresponds to the maximum order of $g^{(i)}(0,1)=0$ in $C_{n,m}$, or the degree of PRM to approach WENO at the endpoint, the larger value of $m$ will favor numerical stability as indicated in Section 3. The choice of $m = n$ is suggested in this study.

(3) The determination of $c_1$

As $c_1$ corresponds to extent of Flatness-II and also *PEC*, a relatively small value is suggested to ensure enough stability and also sufficient space for subsequent regulations. Its recommendation is $c_1 \le 1$ or likewise.

(4) The determination of $c_2$

As $c_2$ corresponds to extent of Flatness-II directly, its value will affect the resolution of mapped WENO scheme. Usually $c_2$ can take the value larger than $10^3 \sim 10^4$ or likewise; besides, the value might differ case by case of schemes and should appear in piecewise manner.

(5) The determination of $n_1$

$n_1=1$ is suggested in consideration of computation efficiency, which favors the larger flatness.

(6) The determination of $m_1$

As $m_1$ is critical to control the transition of PRM from the flat profile around $\omega = d_k$ to identity mapping at endpoint, its proper choice helps to achieve desired flatness and *PEC* concurrently.

Additionally, our practices indicate too large $m_1$ will on the one hand increase the computation cost, on the other hand engender abrupt transition prone to numerical instability.

From above recommendations, the interaction of functionality of $c_1$, $c_2$, and $m_1$ is demonstrated, through which the desired flatness and *PEC* can be achieved possibly. After extensive practices analytically and numerically, the finalized parameters of $\text{PRM}_{n,n}^{L/R,n+1}$ for WENO3,5,7-JS are acquired, which are tabulated in Table 6 corresponding to WENO3,5,7-JS by $r$. Parameters of WENO9-JS are not investigated because of its seldom usage in applications.

Table 6. Parameters of $\text{PRM}_{n,n}^{L/R,n+1}$ by Eq. (16) corresponding to WENO3,5,7-JS with $r$=2-4

|       |       |           |   | $c_1$ | $c_2$ | $m_1$ |
|-------|-------|-----------|---|-------|-------|-------|
| $r=2$ | $n=1$ | $d_0$=1/3 | L | 1 | $7\times10^7$ | 5 |
|       |       |           | R | 1 | $3\times10^6$ | 5 |
|       |       | $d_1$=2/3 | L | 1 | $1\times10^5$ | 4 |
|       |       |           | R | 1 | $3\times10^6$ | 4 |
| $r=3$ | $n=2$ | $d_0$=1/10 | L | 1 | $1\times10^9$ | 5 |
|       |       |           | R | 1 | $5\times10^4$ | 6 |
|       |       | $d_1$=6/10 | L | 1 | $6\times10^5$ | 6 |
|       |       |           | R | 1 | $6\times10^7$ | 6 |
|       |       | $d_2$=3/10 | L | 1 | $3\times10^8$ | 6 |
|       |       |           | R | 1 | $2\times10^5$ | 6 |
| $r=4$ | $n=3$ | $d_0$=1/35 | L | 1 | $1\times10^{11}$ | 5 |
|       |       |           | R | 1 | $5\times10^2$ | 5 |
|       |       | $d_1$=12/35 | L | 1 | $3\times10^4$ | 5 |
|       |       |           | R | 1 | $3\times10^3$ | 4 |
|       |       | $d_2$=18/35 | L | 1 | $1\times10^4$ | 5 |
|       |       |           | R | 1 | $2\times10^4$ | 4 |
|       |       | $d_3$=4/35 | L | 1 | $5\times10^7$ | 5 |
|       |       |           | R | 1 | $5\times10^2$ | 4 |

For illustration, the profiles of integrated $\text{PRM}_{n,n}^{n+1}$ in [0, 1] are shown in Figs. 5-7, where some comparative mappings are drawn as well. Concretely, the profile of $\text{PRM}_{1,1}^2$ for weights of WENO3-JS is displayed in Fig. 5 with the comparison of $\text{PPM}_{1,0}^2$, $g_M$ and $\text{IM}_{2,0;0.1}^3$. Benefited from the strong regulation capability, $\text{PRM}_{1,1}^2$ appears the largest Flatness-II as well as the enhanced *PEC*, which favors the improvement of resolution in smooth region and preservation of stability near discontinuities. Although $\text{IM}_{2,0;0.1}^3$ shows a profile with the least *PEC* which favors the increase of numerical resolution, the distribution might be subjected to potential instability. As shown in Section 5 later, WENO5-$\text{IM}_{2,0;0.1}^3$ cannot pass the computation of blast wave.

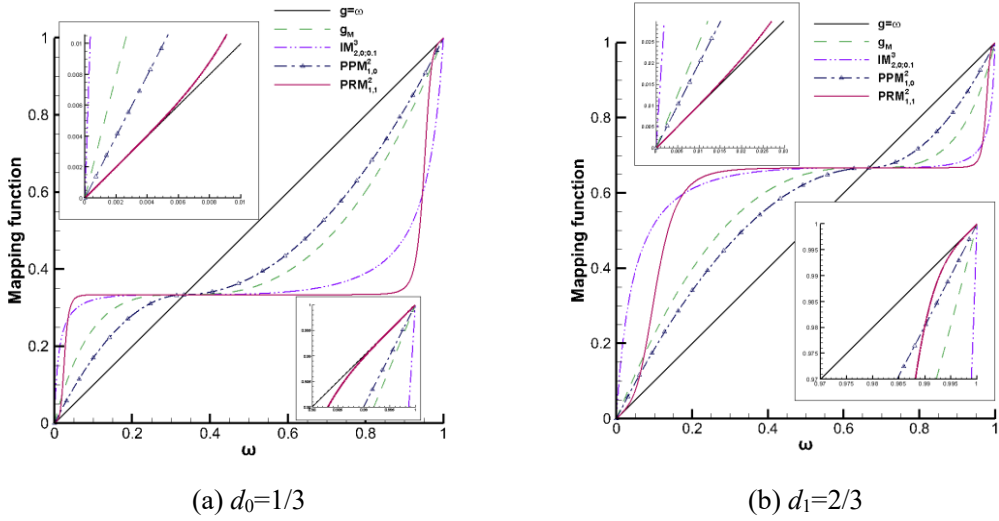

(a) $d_0=1/3$  (b) $d_1=2/3$

Fig. 5. Distributions of integrated $PRM^2_{1,1}$ with comparisons of $PPM^2_{1,0}$, $g_M$ and $IM^3_{2,0;0.1}$; zoomed views regard details near endpoints

Moreover, although $PRM^2_{1,1}$ encompasses the other mappings except near endpoints, we have checked that it appears the least flat in the neighborhood of $\omega=d_k$. As previously mentioned, the latter regards Flatness-I determined by $n$ in $C_{n,m}$. Hence the profiles indicate Flatness-II might differ from Flatness-I in some cases, and the one having larger Flatness-II may have smaller Flatness-I.

In Fig. 8, the profile of integrated $PRM^3_{2,2}$ for weights of WENO5-JS is shown with the comparison of $g_M$, $IM^3_{2,0;0.1}$, $PPM^4_{2,1}$, $R^3_{2,2}$, $PM^8_{6,1}$, and $RM^7_{6,2,0}$, where $R^3_{2,2}$ is one of our practices in point "f" in Section 2 with specific choice of parameters shown in Table 7. The parameters are so chosen that $R^3_{2,2}$ would have the similar Flatness–II as that of $PRM^3_{2,2}$. From the figure, $PRM^3_{2,2}$ is designed to have both the large Flatness-II and high rate of endpoint convergence toward the identity mapping. If comparing with $RM^7_{6,2,0}$, one may find that $PRM^3_{2,2}$ on the one hand have the larger Flatness-II, on the other hand possesses the especially enforced *PEC* at endpoints. Such regulation can only be achieved by means of Eqns. (14) -(16) with the proper choices of parameters. Particularly, one can see that in the case of $d_1=6/10$, the profile of $PRM^3_{2,2}$ reflects a well balance of flatness and *PEC*, while $RM^7_{6,2,0}$ exhibits an exaggerated flat profile and slower rate of endpoint convergence near the endpoint $\omega = 1$. Hence the advantage of profile from the perspective of mathematics of $PRM^3_{2,2}$ is demonstrated. Regarding $R^3_{2,2}$, although it has the similar Flatness-II as that of $PRM^3_{2,2}$, the mapping indicates *PEC* relatively inferior. As shown in Section 5.3, $R^3_{2,2}$ will yield oscillations in the long time advection of combination wave problem, which validates the importance of *PEC* on numerical stability.

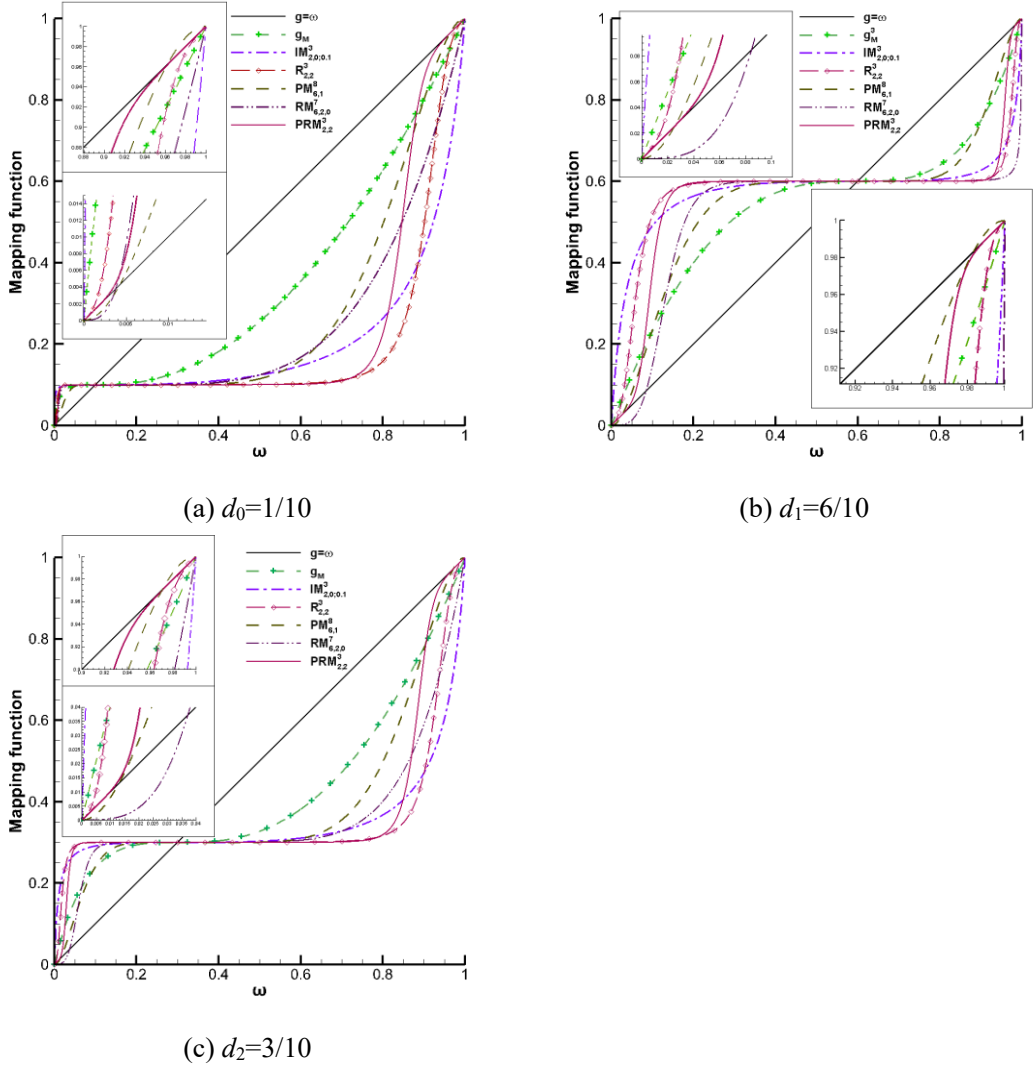

(a) $d_0=1/10$

(b) $d_1=6/10$

(c) $d_2=3/10$

Fig. 6. Distributions of integrated $PRM^3_{2,2}$ with $g_M$, $IM^3_{2,0;0.1}$, $R^3_{2,2}$, $PM^8_{6,1}$, and $RM^7_{6,2,0}$ with zoomed views regarding details near endpoints

The profile of integrated $PRM^4_{3,3}$ for weights of WENO7-JS is shown in Fig. 7 with the comparisons of $g_M$, $IM^3_{2,0;0.1}$, $PM^8_{6,1}$, and $RM^7_{6,2,0}$. With the help of sufficient capability of regulation, $PRM^4_{3,3}$ is designed to have a profile with the balanced large Flatness-II and well-performed *PEC* which outweighs the comparatives.

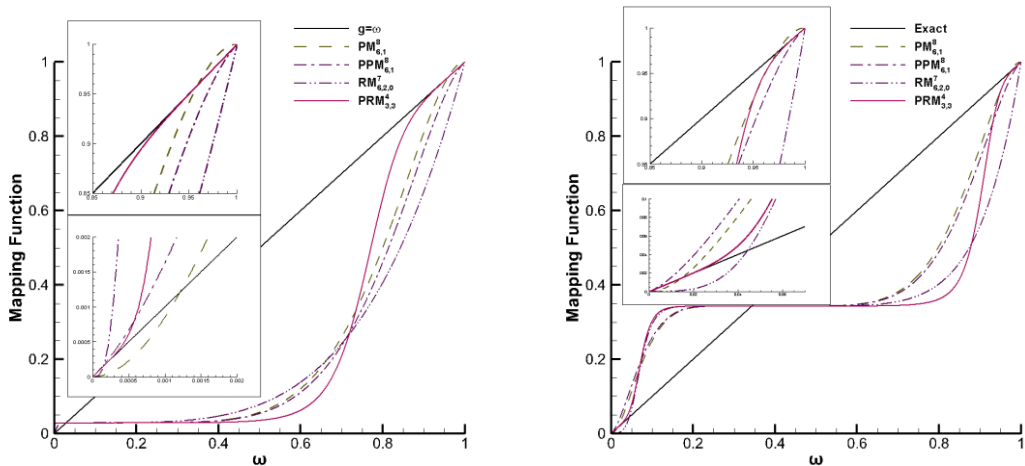

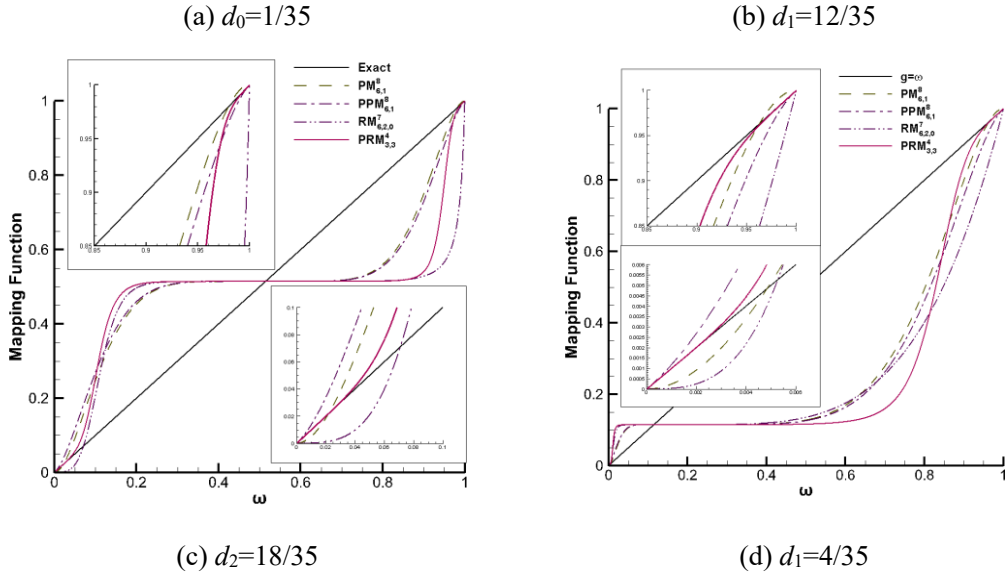

(a) $d_0=1/35$                (b) $d_1=12/35$

(c) $d_2=18/35$             (d) $d_1=4/35$

Fig. 7. Distributions of integrated $PRM^3_{2,2}$ with $g_M$, $PM^8_{6,1}$, $PPM^8_{6,1}$ and $RM^7_{6,2,0}$ with zoomed views regarding details near endpoints

Finally, in order to demonstrate the sufficient regulation capability of PRM, a particular $PRM^3_{2,2}$, namely $PRM^{3*}_{2,2}$, is provided as an example to mimic profiles of $PPM^8_{6,1}$ and $RM^7_{6,2,0}$ in the case of $d_k=6/10$ through specific choice of parameters. Corresponding results are shown in Fig. 8, and the parameters regarding $PRM^3_{2,2}$ are also given in Table 7.

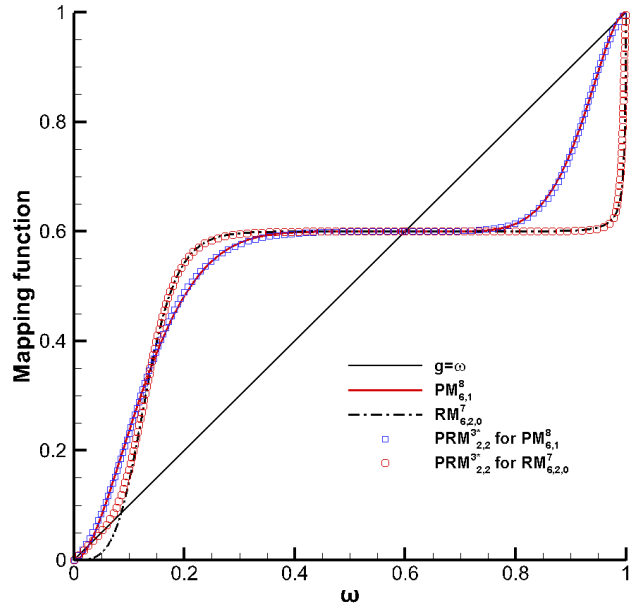

Fig. 8. Imitation of $PM^8_{6,1}$ and $RM^7_{6,2,0}$ by $PRM^3_{2,2}$ in the case of $d_k=6/10$ through specific choice of parameters

Table 7. Parameters of $R^3_{2,2}$ that resembles Flatness-II of $PRM^3_{2,2}$ and that of $PRM^{3*}_{2,2}$ to mimic $PM^8_{6,1}$ and $RM^7_{6,2,0}$ in the case of $d_k=6/10$ with mapping defined by Eq. (16)

|  |  |  | $c_1$ | $c_2$ | $m_1$ |
|---|---|---|---|---|---|
| $R^3_{2,2}$ | $d_0=1/10$ | L | 30090 | 0 | 0 |
|  |  | R | 676.6666 | 0 | 0 |

|  | $d_1$=6/10 | L | 1235.6790 | 0 | 0 |
|---|---|---|---|---|---|
|  |  | R | 8335 | 0 | 0 |
|  | $d_2$=3/10 | L | 12970.7047 | 0 | 0 |
|  |  | R | 929.2592 | 0 | 0 |
| PRM$^{3*}_{2,2}$ ⇒ PM$^8_{6,1}$ | $d_1$=6/10 | L | 26 | 13 | 2 |
|  |  | R | 40 | 20 | 2 |
| PRM$^{3*}_{2,2}$ ⇒ RM$^7_{6,2,0}$ | $d_1$=6/10 | L | 1 | 7500 | 5 |
|  |  | R | 1000 | 10000 | 2 |

## 5 Numerical tests

5.1 Case descriptions

The following 1-D problems of scalar advection and that by Euler equations are chosen to explore the performance of PRMs.

(1) 1-D scalar advection problems

The governing equation is: $\partial u/\partial t + \partial u/\partial x = 0$ with various conditions $u(x,0)$ representing specific problems. Concretely, the initial conditions are:

(a) Sinusoidal-like wave advection I (SWA-I): $u(x,0) = \sin\left(\pi x - \frac{\sin(\pi x)}{a\pi}\right)$ at $x \in [-1,1]$

The case describes advection of sinusoidal-like wave with only two first-order critical points when $a > 1/\pi$, and the choice of $a=1$ has been used in Ref. 3. Other than canonical one as $u(x,0) = \sin(\pi x)$, we find SWA-I would make numerical schemes behave in quite disparate convergence rates of accuracy. Empirically, the scheme having fast convergence rate is expected to achieve desired accuracy in less grids, which is rather favored by applications. Specifically, $a=1$ and $1.005/\pi$ are chosen for the third and fifth-order scheme case respectively in Section 5.3. The fourth-order Runge-Kutta method is employed for temporal discretization with $\Delta t < \Delta x^{\frac{2r-1}{4}}$. A series of grids with the numbers {20, 40, 80, …} are employed and sequentially referred as the first, second, … grids for convenience. The evolution of L$_\infty$-norm errors with $\Delta x$ after one period of advection is acquired to indicate the accuracy order.

(b) Sinusoidal-like wave advection II (SWA-II): $u(x,0) = \sin^3(\pi x - \frac{\sin(\pi x)}{a\pi})$ at $x \in [-1,1]$

This case is also designed to test the order accuracy of schemes. Besides two first-order critical points, this case has two second-order critical points at $x=0, \pm1$. $a=0.32$ is chosen herein with the first-order critical points at $x \approx \pm 0.7345$. Numerical schemes will usually show more disparate rates of order-convergence than in the case of $u(x,0) = \sin^3(\pi x)$. The temporal scheme and computation configurations such as grids, $\Delta t$ and computation period are the same as that in SWA-I.

(c) Combination-waves advection [4]

$$u(x,0) = \begin{cases} \frac{1}{6}(G(x,\beta,z-\delta)+G(x,\beta,z+\delta)+4G(x,\beta,z)), & -0.8 \le x \le -0.6 \\ 1, & -0.4 \le x \le -0.2 \\ 1-|10(x-0.1)|, & 0 \le x \le 0.2 \\ \frac{1}{6}(F(x,\alpha,a-\delta)+F(x,\alpha,a+\delta)+4F(x,\alpha,a)), & 0.4 \le x \le 0.6 \\ 0, & \text{otherwise} \end{cases} \quad \text{at} \quad x \in [-1,1],$$

where $G(x,\beta,z) = e^{-\beta(x-z)^2}$, $F(x,\alpha,a) = \sqrt{\max(1-\alpha^2(x-a)^2, 0)}$, $a = 0.5$, $z = -0.7$, $\delta = 0.005$, $\alpha = 10$ and $\beta = \log 2/36\delta^2$. Typical grid numbers include: 200, 400 and 800 as in Refs. 7-8, 10-11. Computational time can vary from the short period 2 to the long one 4000. In the computations, the third-order TVD Runge-Kutta method [2] is employed with CFL number equals 0.1.

This example has already served as canonical test to check numerical stability in short time period, however the stability as well as corresponding errors of long time computation was first concerned by Ref. 7 regarding WENO5 schemes. Later, computations in short and long periods were intensively investigated in Refs. 8, 10 and 11 regarding mapped WENO5 and WENO7. Because of the importance of referred issues in practical applications, they are studied in detail in this study.

(2) 1-D problems by Euler equations

Because all schemes discussed next can pass standard tests such as Sod problem, such practices will not be revisited here for brevity. The problems are chosen as follows:

(a) Strong shock wave [18]

This case poses a trial regarding the robustness of schemes. The initial condition is:

$$(\rho, u, p) = \begin{cases} (1, 0, 0.1PR) & -5 \le x < 0 \\ (1, 0, 0.1) & 0 < x \le 5 \end{cases}.$$ In this study, $PR$ takes the value of $10^3$ and $10^6$ respectively, and corresponding computations advance to $t=0.3$ and $0.01$. Typical grid number is $N=201$. Because the schemes that fulfill the computation with $PR=10^3$ are found to yield quite similar distributions, only results of $PR=10^6$ will be shown in subsequent studies.

(b) Blast wave

This canonical example raises as another trial regarding the robustness of schemes. The initial condition is: $(\rho, u, p) = \begin{cases} (1, 0, 1000) & 0 \le x < 0.1 \\ (1, 0, 0.01) & 0.1 \le x \le 0.9 \\ (1, 0, 100) & 0.9 < x \le 1 \end{cases}$, and solid wall condition is imposed at boundaries. The grid number is $N=200$ and the computation advances to $t=0.038$. By convention, a result at 100001 grids by WENO5-JS is regarded as the "Exact" solution for reference.

(c) Shu-Osher problem

This problem is a benchmark to test the resolution of numerical schemes. The initial condition is: $(\rho, u, p) = \begin{cases} (3.857143, 2.629369, 10.3333) & -5 \le x < -4 \\ (1+0.2\sin(5x), 0, 1) & -4 < x \le 5 \end{cases}$. The computation advances to $t=1.8$. By convention, the result of WENO5-JS at 2001 grids is regarded as the "Exact" solution for reference.

(d) Titarev-Toro problem

This example describes a Mach 1.1 shock interacting with the density fluctuations with high frequency, which serves as another trial to test the resolution of numerical schemes. The initial condition is: $(\rho, u, p) = \begin{cases} (1.515695, 0.523346, 1.805) & -5 \leq x < -4.5 \\ (1+0.1\sin(20\pi x), 0, 1) & -4.5 \leq x < 5 \end{cases}$. Typical grid number is $N=1000$ and the computation runs to $t=5$. By convention, the result of WENO5-JS at 10001 grids is referenced as the "Exact" solution.

5.2 Additional WENO-Z type schemes for comparison

In order to further explore the performance of *PRM* methods, not only mappings in Section 2.2 are compared but also comparisons are made with some typical/updated WENO-Z type schemes. Under the framework of Eqs. (3) -(5), their formulations are described next for completeness.

(1) Third-order schemes: WENO3-P+3 [14] and WENO3-F3 [15]

WENO3-P+3 is obtained by applying new $\alpha_k$ and $\tau_p$ as $\alpha_k = d_k \left[ 1 + \dfrac{\tau_p}{IS_k^{(2)} + \varepsilon} + \lambda \left( \dfrac{IS_k^{(2)} + \varepsilon}{\tau_p + \varepsilon} \right) \right]$ with $\tau_p = \left| \tfrac{1}{2}(IS_0^{(2)} + IS_1^{(2)}) - \tfrac{1}{4}(f_{i-1} - f_{i+1})^2 \right|$, $\lambda = (\Delta x)^{1/6}$ and $\varepsilon = 10^{-40}$.

WENO3-F3 is acquired by applying new $\alpha_k$, $IS_3$ and $\tau_{F3}$ as $\alpha_k = d_k \left( 1 + \dfrac{\tau_{F3}}{IS_k^{(2)} + \varepsilon} \right)$ where

$\tau_{F3} = \left| \tfrac{1}{2}(IS_0^{(2)} + IS_1^{(2)}) - IS_3 \right|^p$, $IS_3 = \tfrac{1}{12}(f_{i-1} - 2f_i + f_{i+1})^2 + \tfrac{1}{4}(f_{i-1} - f_{i+1})^2$, $p=3/2$ and $\varepsilon = 10^{-40}$.

(2) Fifth-order schemes: canonical WENO5-Z scheme [4] and WENO5-NIS [16]

Canonical WENO5-Z scheme is derived by applying new $\alpha_k$ and $\tau_z$ as $\alpha_k = d_k \left[ 1 + \left( \dfrac{\tau_z}{IS_k^{(3)} + \varepsilon} \right)^q \right]$ where $\tau_z = \left| IS_0^{(3)} - IS_2^{(3)} \right|$. When the first-order critical points occur, WENO5-Z would have the fourth-order accuracy with $q=1$ but preserve the fifth-order with $q=2$, where the choice of $q=1$ is usually thought to yield higher resolution [4].

WENO5-NIS is formulated by proposing the new smoothness indicators $NIS_k^{(3)}$ from the original $IS_k^{(3)}$ as $\begin{pmatrix} NIS_0^{(3)} \\ NIS_1^{(3)} \\ NIS_2^{(3)} \end{pmatrix} = \begin{pmatrix} IS_0^{(3)} - \left|(f_i - 2f_{i+1} + f_{i+2})(3f_i - 4f_{i+1} + f_{i+2})\right| \\ IS_1^{(3)} - \left|(f_{i-1} - 2f_i + f_{i+1})(f_{i-1} - f_{i+1})\right| \\ IS_2^{(3)} - \left|(f_{i-2} - 2f_{i-1} + f_i)(f_{i-2} - 4f_{i-1} + 3f_i)\right| \end{pmatrix}$.

5.3 1-D scalar problems

(1) Results of WENO3-PRM$^2_{1,1}$ and corresponding third-order comparatives

(a) Accuracy preservation at critical points and convergence rate by SWA-I

In this situation, SWA-I is tested by WENO3-PRM$^2_{1,1}$ at the choice of $a=1$ where the first-order critical points exist. For comparison, canonical WENO3-JS and other mapped WENOs are tested as well. It is worth mentioning that the smoothness indicators in the mappings, including that in PRM$^2_{1,1}$, should employ the ones in WENO5-JS as described in Section 3. Comparative mappings include: PPM$^2_{1,0}$, $g_M$, PPM$^3_{2,0}$ and IM$^3_{2,0.1}$. PPM$^2_{1,0}$ has the minimum requirement of $r_{c,g}=1$ as that of PRM$^2_{1,1}$, whereas the rest have one order higher as $r_{c,g}=2$.

Figure 9 illustrates that all mapped schemes can recover the third-order accuracy, however, different schemes demonstrate disparate rates of convergence. Among them, WENO3-*PRM*$^2_{1,1}$ shows the quickest rate which achieves the 2.8$^{th}$-order at the second grids and third-order at the third grids, in the meanwhile WENO3-IM$^3_{2,0.1}$ performs similarly but with slightly smaller accuracy order at first two grids. The other mappings show definite slower rate of convergence even if $g_M$ and PPM$^3_{2,0}$ have one order larger of $r_{c,g}$ than that of PRM$^2_{1,1}$. As expected, WENO3-JS shows a degraded order of 2 because of the critical points.

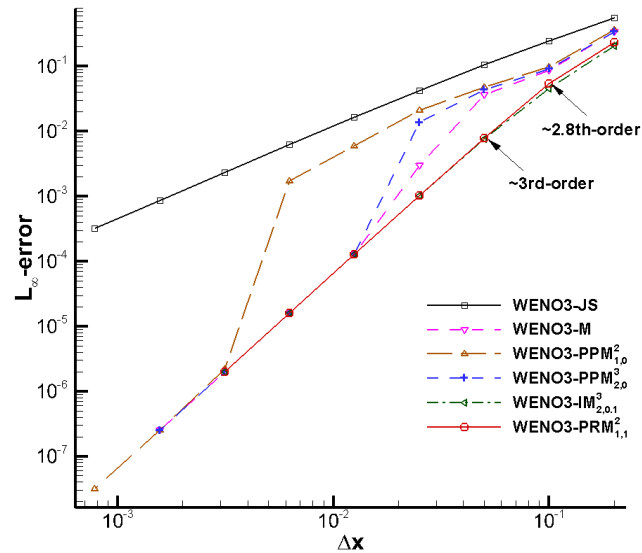

Fig. 9. Evolution of accuracy convergence in the case of scalar advection with the intial condition $u(x,0) = \sin\left(\pi x - \frac{\sin(\pi x)}{\pi}\right)$ by WENO3-JS, -M, -PPM$^2_{1,0}$, -PPM$^3_{2,0}$, -IM$^3_{2,0.1}$ and –PRM$^2_{1,1}$

(b) Long time computation of combination-waves advection

A long period up to T=4000 is adopted in computation. For comparison, the following mappings are selected: $g_M$, PPM$^3_{2,0}$ and IM$^3_{2,0.1}$, whose $r_{c,g}$ order is 2 and one larger than that of PRM$^2_{1,1}$; likewise, WENO3-JS is also included. In Fig. 10, the distribution of WENO3-PRM$^2_{1,1}$ on 800 grids is drawn with the comparisons by other schemes. One can see that after long time computation, WENO3-PRM$^2_{1,1}$ and -IM$^3_{2,0;0.1}$ show quite similar results, which are relatively more accurate than that by WENO3-M and PPM$^3_{2,0}$, whereas WENO3-JS (with $\varepsilon=10^{-6}$ in Eq. (4)) yield a rather disspative distribution. From the figure, no instability occurs and all schemes perform stably.

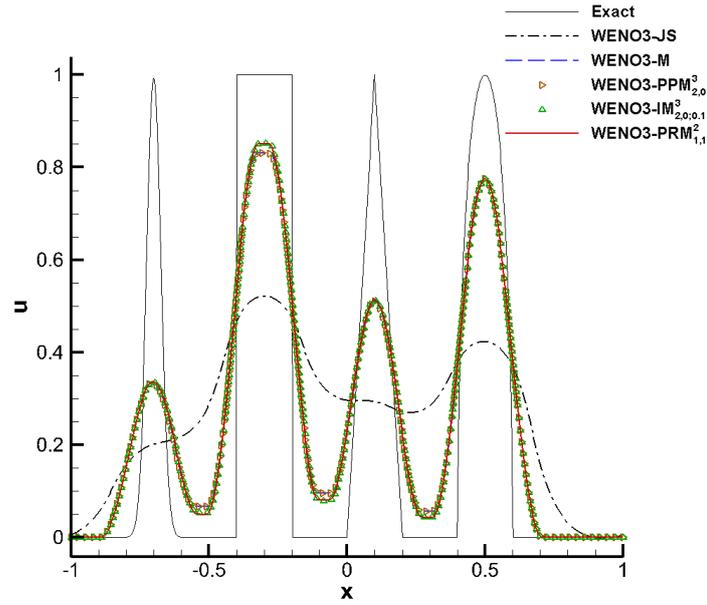

Fig. 10. Results of combination-waves advection by WENO3-PRM$^2_{1,1}$ with the comparisons by WENO3-JS, M, -PPM $^3_{2,0}$ and -IM $^3_{2,0.1}$ on 800 grids and with $T$=4000

In order to reveal numerical errors of schemes, additional 200, 400 grids are adopted and corresponding $L_1$-error is computed. As shown in Fig. 11, WENO3-PRM$^2_{1,1}$ indicates the least error among four mapped WENOs, whereas WENO3-JS exhibits an obvious large error. Moreover, carefully checking shows WENO3-PRM$^2_{1,1}$ achieves the approximate first-order on the third grids.

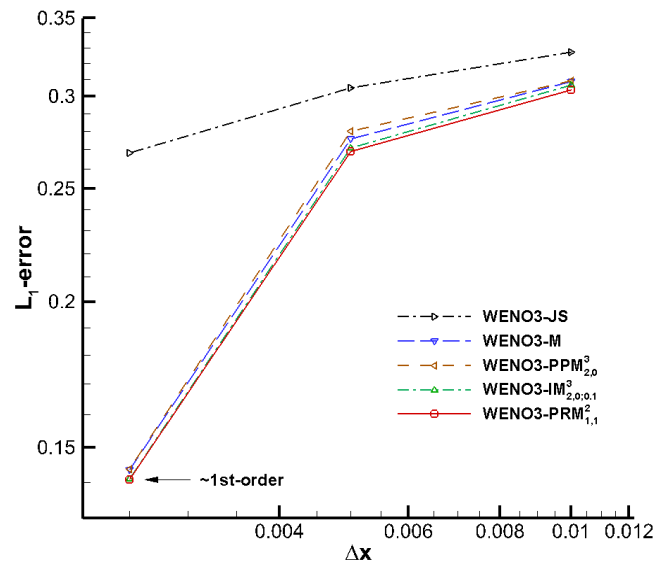

Fig. 11. $L_1$-errors with $\Delta x$ from computations of combination-waves advection at $T$=4000 by WENO3-PRM $^2_{1,1}$ with the comparisons by WENO3-M, -PPM $^3_{2,0}$ and -IM $^3_{2,0.1}$

(2) Results of WENO5-PRM$^3_{2,2}$ and corresponding fifth-order comparatives
(a) Accuracy preservation at critical points and convergence rate by SWA-I

As shown in Section 3 and Table 12 in Appendix I, the minimum requirement of $r_{c,g}$ to preserve the optimal order is 2 and the largest $n_{cp}$ where 5$^{\text{th}}$-order be preserved is 1 for mapped WENO5-JS. Therefore, WENO5-PRM$^3_{2,2}$ is still tested by SWA-I to check accuracy order and convergence rate, where $a$ is taken as $1.005/\pi$ in the initial condition. For comparison, the following mappings with the same $r_{c,g}$ =2 are tested, i.e. $g_M$, PPM$^3_{2,0}$, and IM$^3_{2,0.1}$; besides, mappings with larger $r_{c,g}$ proposed in Refs. 8 and 10 for WENO5-JS are also checked, namely PM$^8_{6,1}$ and RM$^7_{6,2,0}$. It is worthwhile to mention that AIM$^5_{4,2;1E4}$ fails to work in the case of WENO5-JS probably because of its parameters, namely $c=10^4$, only being applicable to WENO7-JS. According to analysis above and results in Fig. 12, all mapped schemes achieve the fifth-order theoretically and numerically. Furthermore, PRM$^3_{2,2}$ demonstrates the quickest convergence rate which reaches the optimal order on the third grids, whereas the others achieve the order at the fifth or sixth grids ($g_M$). As expected, WENO5-JS can only obtain a degraded third order from the figure.

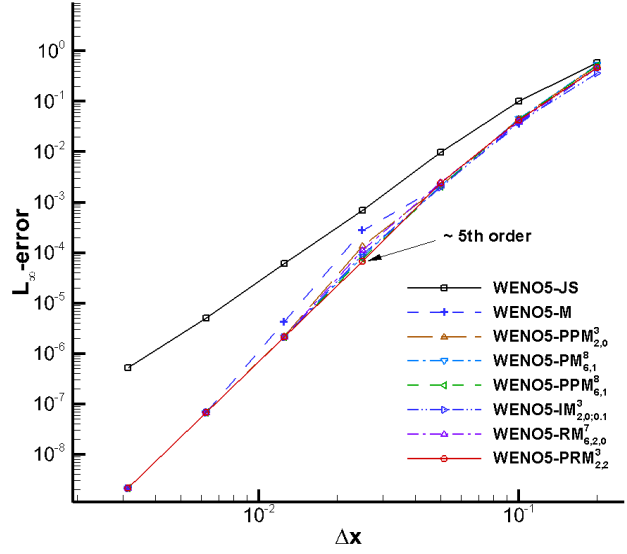

Fig. 12. Evolution of accuracy convergence in the case of scalar advection with the intial condition $u(x,0) = \sin\left(\pi x - \frac{\sin(\pi x)}{1.005}\right)$ by WENO5-JS, -M, -PPM$^3_{2,0}$, -IM$^3_{2,0;0.1}$, -PM$^8_{6,1}$, -RM$^7_{6,2,0}$ and –PRM$^3_{2,2}$

(b) Long time computation of combination-waves advection

WENO5-PRM$^3_{2,2}$ is tested on 200, 400 and 800 grids with $T$=2000. In Refs. 7, 8 and 10, the stability and numerical errors of this case were intensively studied, where WENO5-RM$^7_{6,2,0}$ was considered to have outstanding performances both in short time and long time compuations. For comparison, the following schemes are checked: WENO5-JS, -M, -IM$^3_{2,0;0.1}$, -PM$^8_{6,1}$ and -RM$^7_{6,2,0}$, where AIM $^5_{4,2;1E4}$ is absent because of its failure in computation. Corresponding results are shown in Fig. 13, and on checking, the results of schemes for comparison coninicde with that in Ref. 10. The figure tells that WENO-JS and –M show smeared distributions after long time computation, WENO5-IM$^3_{2,0;0.1}$ yields relatively disspiative result (see the distribution at $x$=-0.7 and -0.8) and the rest schemes perform similarly except at the corner $x$=0.6. At the corner, WENO5-IM$^3_{2,0;0.1}$ and -RM$^7_{6,2,0}$ yield results without oscillations, whereas WENO5-PM$^8_{6,1}$ and -PRM$^3_{2,2}$ show smaller ones. As implied Refs. 3 and 10, $\varepsilon$ in eq. (4) may affect the stability of computation especially that in the long time run. As shown in Fig. 13(c), one can see that WENO5-RM$^7_{6,2,0}$ yields similar oscillations with $\varepsilon=10^{-45}$ and that of WENO5-PRM$^3_{2,2}$ can be totally free of fluctuations with $\varepsilon=10^{-101}$.

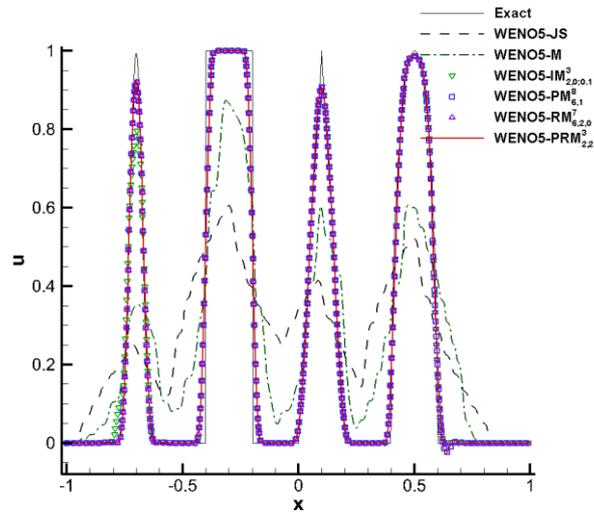

(a) Global view

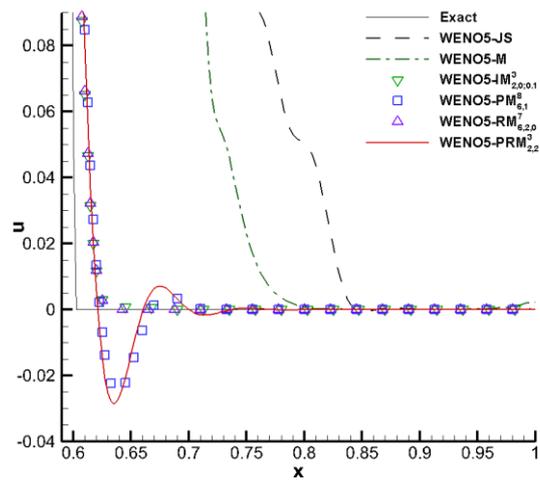

(b) Zoomed view

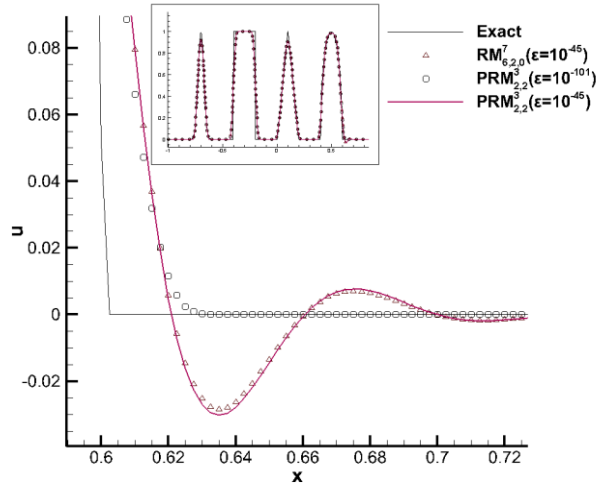

(c) Results of PRM$^3_{2,2}$ with different $\varepsilon$ in Eq. (4) and RM$^7_{6,2,0}$ with $\varepsilon=10^{-45}$

Fig. 13. Results of combination-waves advection I by WENO5-PRM$^3_{2,2}$ with the comparisons by that of WENO5-JS, -M, -IM$^3_{2,0;0.1}$, -PM$^8_{6,1}$ and -RM$^7_{6,2,0}$ on 800 grids and with $T=2000$

Next, L$_1$-error of WENO5- PRM$^3_{2,2}$ is invesitgated with the comparisons by WENO5-JS, -M, -IM$^3_{2,0;0.1}$, -PM$^8_{6,1}$ and -RM$^7_{6,2,0}$. As shown in Fig. 14, WENO5-PRM$^3_{2,2}$ indicates a low error level as that of WENO5-PM$^8_{6,1}$ and -RM$^7_{6,2,0}$, and they nearly achieve the first-order on the second grids. As the comparison, WENO5-IM$^3_{2,0;0.1}$ shows a relatively larger error whereas that of WENO-M and WENO5-JS fall into the group with the largest errors.

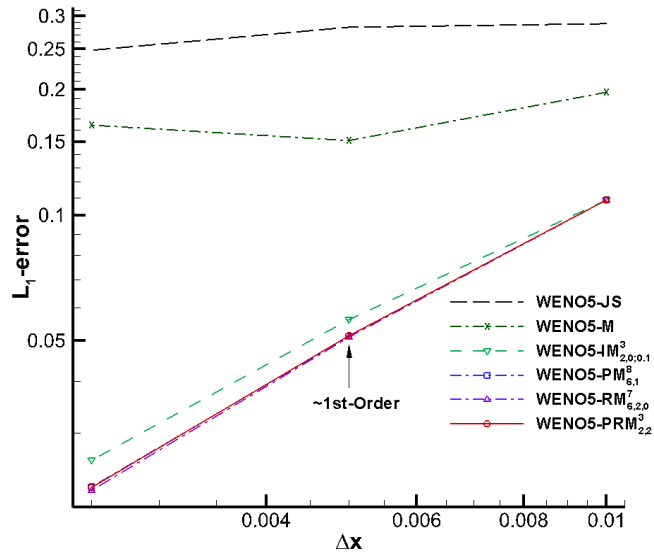

Fig. 14. L$_1$-errors with $\Delta x$ from computations of combination-waves advection at $T=2000$ by WENO5-PRM$^3_{2,2}$ with the comparisons by WENO5-JS, -M, -IM$^3_{2;0.1}$, -PM$^8_{6,1}$ and -RM$^7_{6,2,0}$

In Section 4, a comparative R$^3_{2,2}$ is devised which has the similar Flatness–II as that of PRM$^3_{2,2}$. As shown in Fig. 6, R$^3_{2,2}$ inidcates relatively inferior performances on *PEC* due to the lack of sufficient regulation. Two mappings are also used in current computation at $T=2000$ on 800 grids and under $\varepsilon=10^{-45}$ in Eq. (4), and the results is shown in Fig. 15. The figure tells that R$^3_{2,2}$ yields distinct oscillations than PRM$^3_{2,2}$ does, which validates the functionality of new implementation in PRM$^3_{2,2}$ in numerical stability enhancement and verify the importance of *PEC*.

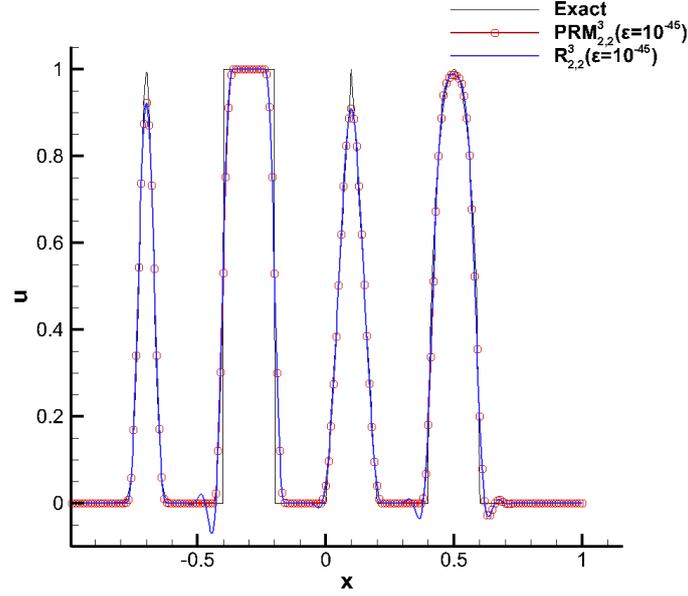

Fig. 16. Comparisons of WENO5-PRM$^3_{2,2}$ with -R$^3_{2,2}$ in combination-waves advection at *T*=2000 on 800 grids and under $\varepsilon=10^{-45}$ in Eq. (4)

(3) Results of WENO7-PRM$^4_{3,3}$ and corresponding seventh-order comparatives

(a) Accuracy preservation at critical points and convergence rate by SWA-II

Recalling the analysis in Section 3 and Table 12 in Appendix I, SWA-II is chosen with the existence of two groups of first- and two second-order critical points. The case represents an most possible situation that canonical WENO7-JS could preserve its optimal order by mappings, and the minimum requirement of $r_{c,g}$ for order preservation is 3. To compare with WENO7-PRM$^4_{3,3}$, canonical WENO7-JS and some recent mappings are chosen: PM$^8_{6,1}$, RM$^7_{6,2,0}$ and AIM$^4_{3,3;1E4}$. The reason to choose three mappings is that they are specially designed and/or tested for WENO7-JS[10, 11], while aforementioned IM$^3_{2,0;0.1}$ is not chosen because it does not fulfill the required thrid order of $r_{c,g}$. The evolution of L$_1$-errors with $\Delta x$ are shown in Fig. 16, which indicates all schemes have achieved the optimal orders in the presence of critical points. Carefully checking shows WENO7-PM$^8_{6,1}$ acquires the seventh-order on the fourth grids, while the rest do on the fifth grids. Although WENO7-PM$^8_{6,1}$ exhibits relatively faster rate of order convergence, it indicates the largest error on the third grids.

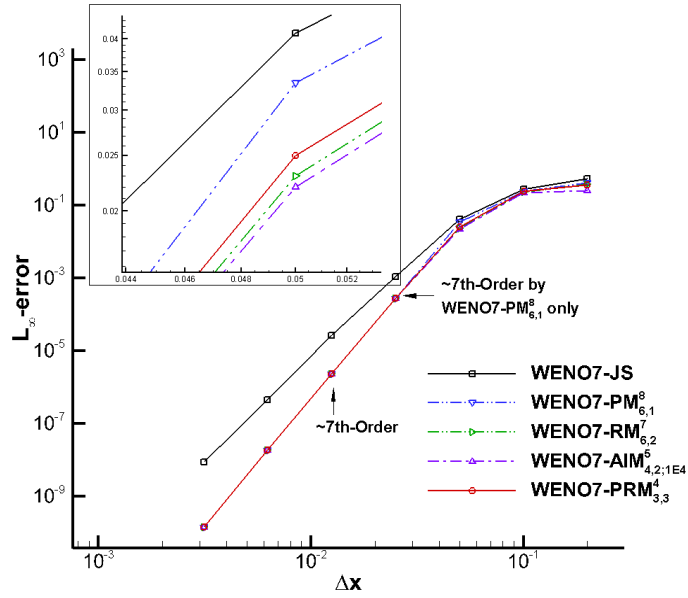

Fig. 16. Evolution of accuracy convergence in the case of scalar advection with the intial condition $u(x,0)=\sin^3(\pi x - \frac{\sin(\pi x)}{0.32\pi})$ by WENO7-PRM$^4_{3,3}$ with the comparisons by WENO7-JS, -PM$^8_{6,1}$, -RM$^7_{6,2,0}$ and -AIM$^5_{4,2;1E4}$.

(b) Long time computation of combination-waves advection

This case is especially concerned and intensively tested among the seventh-order mapped WENOs in Refs. 10-11. For comprehensiveness, the longest computation time from the references is chosen, namely $T=2000$, and the full series of grids therein are used with numbers {200, 400, 800}. Besides WENO7-PRM$^4_{3,3}$, similar schemes are chosen as the previous example. The results of 200 grids is first shown in Fig. 17. Prior to further discussion, it is worthy of mention that we have repeated and verified the result of WENO7-PM$^8_{6,1}$ on 200 grids and $T=1000$ in Ref. 10. From Fig. 17, all schemes have shown discrepencies with respect to the exact solution, where WENO7-AIM$^5_{4,2;1E4}$ yields an obvious under-estimation of the second square and overshoots about the fourth oval, and the other schemes respectively yield small overshoots and obvsious under-estimations also. Because the issue of numerical stablity is the focal in Ref. 10-11, the distributions at the corner $x\approx-0.78$ is zoomed in Fig. 21(b). The figure tells that both WENO7-PM$^8_{6,1}$ and -RM$^7_{6,2,0}$ show osccilations while WENO7-AIM$^5_{4,2;1E4}$ is free of them; besides, WENO7-PRM$^4_{3,3}$ indicates large deviation from the exact solution there. Such deviation would definitely increase computation error, however they should be attributed to numerical dissipation other than instability according to our experience and discussions in Ref. 11. In short, in the case of 200 grids, WENO7-PRM$^4_{3,3}$ and AIM$^5_{4,2;1E4}$ yield results without osillations at T=2000.

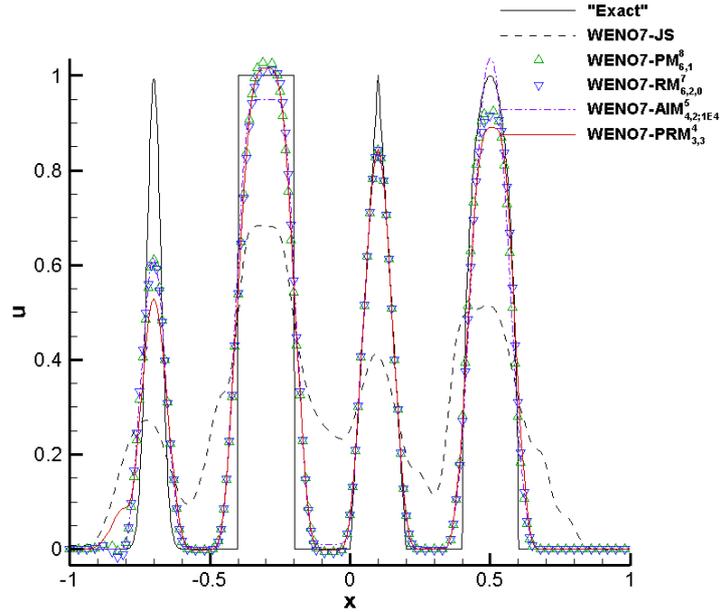

(a) Global view

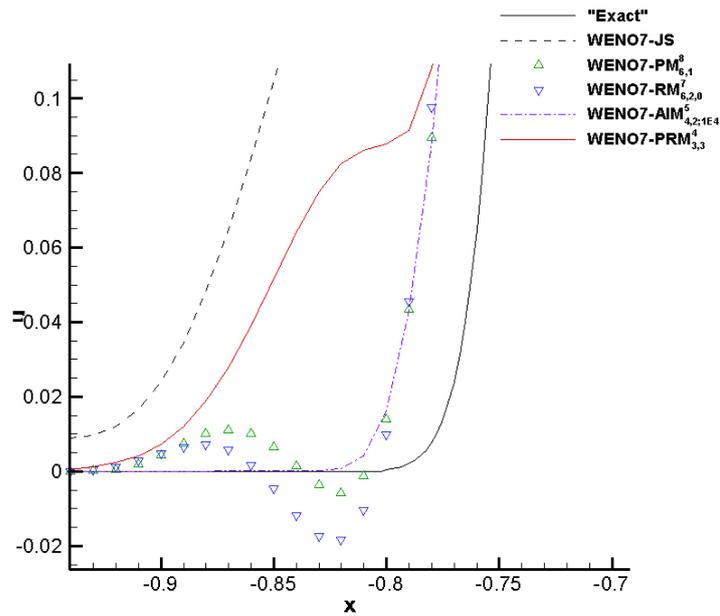

(b) Zoomed view at $\in[-0.96, -0.69]$

Fig. 17. Results of combination-waves advection on 200 grids at $T=2000$ by WENO7-PRM$^4_{3,3}$ with the comparisons by WENO7-JS, -PM$^8_{6,1}$, -RM$^7_{6,2,0}$ and -AIM$^5_{4,2;1E4}$.

Next, the computations are carried out on 400 grids with $T=2000$, and the results are shown in Fig. 18. All mapped schemes preform normally except at the foot of the first peak and the region between the first two structures where oscillations occur. Besides, WENO7-PRM$^4_{3,3}$ yields a relatively smeared description on the first peak. In order to clearly visualize the oscillations, the zoomed view is shown in Fig. 18(b). The figure tells that at the left foot of the first peak, WENO7-PM$^8_{6,1}$, -RM$^7_{6,2,0}$ and -AIM$^5_{4,2;1E4}$ generate oscillations whereas WENO7-PRM$^4_{3,3}$ presents smooth distribution; at the right foot of the peak, all mapped schemes yield oscillations except WENO7-AIM$^5_{4,2;1E4}$; additionally, WENO7-PRM$^4_{3,3}$ produces more perturbations between the first two structures. Overall, all schemes yield oscillations in this case, whereas WENO7-PRM$^4_{3,3}$ might indicate more error because of its deviation from the exact solution at the first peak.

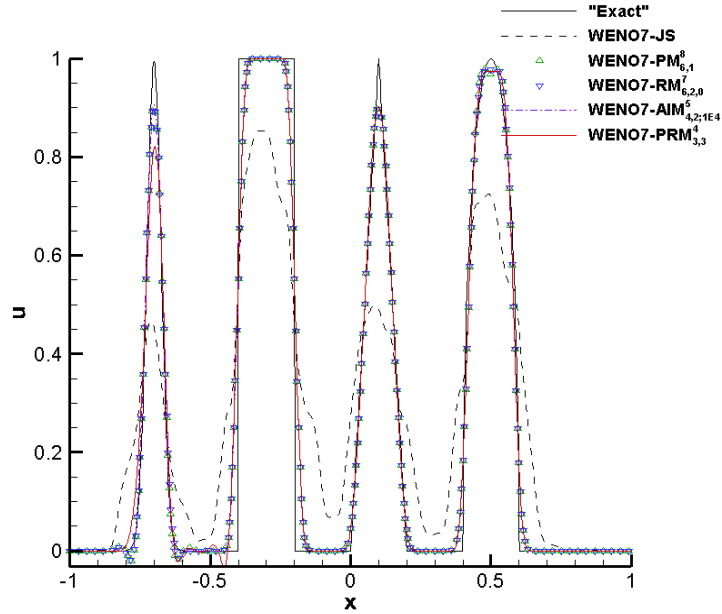

(a) Global view

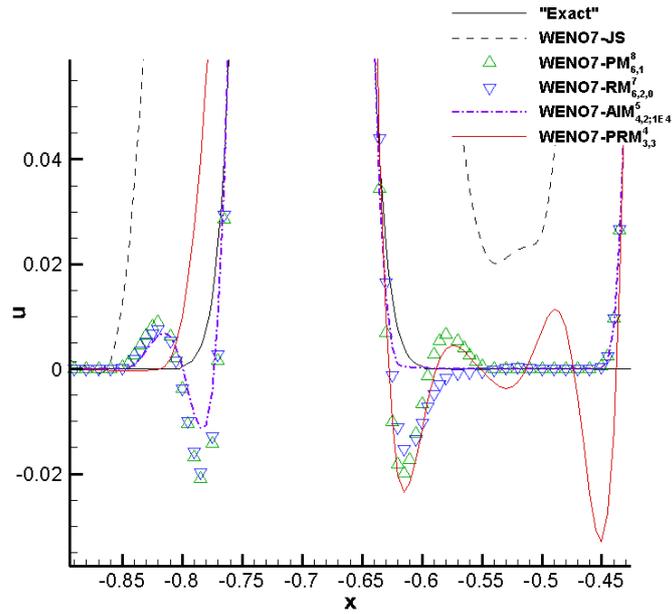

(b) Zoomed view where $x \in [-0.89, -0.43]$

Fig. 18. Results of combination-waves advection on 400 grids at $T=2000$ by WENO7-PRM$^4_{3,3}$ with the comparisons by WENO7-JS, -PM$^8_{6,1}$, -RM$^7_{6,2,0}$ and -AIM$^5_{4,2;1E4}$.

Thirdly, the schemes are tested on 800 grids with $T=2000$. As shown in Fig. 19, all mapped schemes perform nicely except small perturbations by some schemes appear at the right foot of the fourth oval. The zoomed view in Fig. 19(a) indicates that WENO7-PRM$^4_{3,3}$, -PM$^8_{6,1}$ and -RM$^7_{6,2,0}$ take this behavior, whereas WENO7-AIM$^5_{4,2;1E4}$ yields a distribution free of oscillations coinciding with reports in Ref. 11. It is worth mentioning that in our computation, WENO7-RM$^7_{6,2,0}$ performs well in the rest region, and the oscillations at the foots of the second and fourth structures reported in Ref. 11 do not happen.

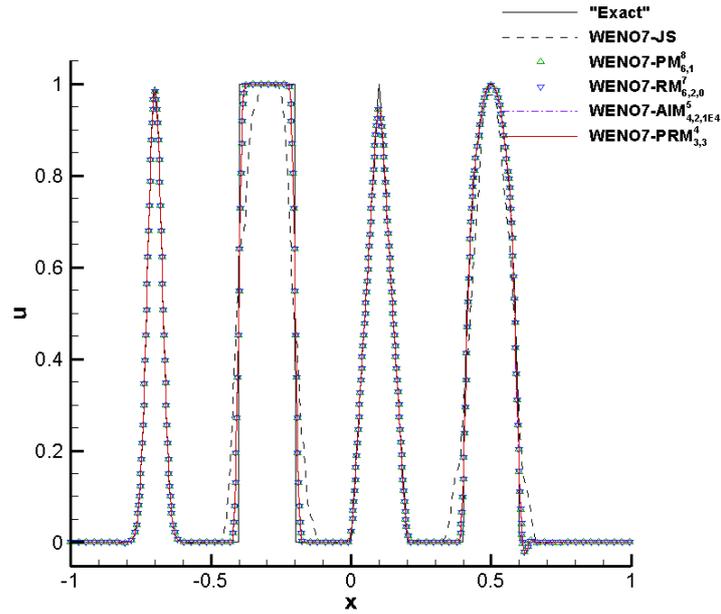

(a) Global view

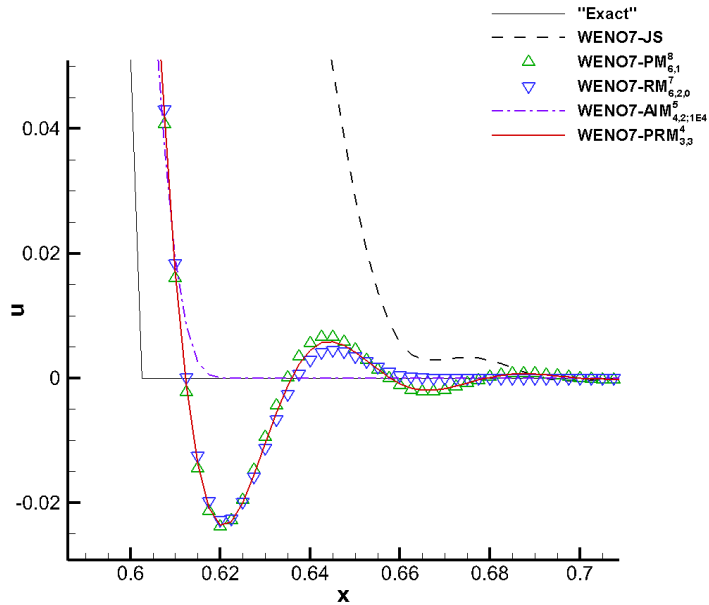

(b) Zoomed view where $x \in [0.59, 0.71]$

Fig. 19. Results of combination-waves advection on 800 grids and at $T=2000$ by WENO7-PRM$^4_{3,3}$ with the comparisons by WENO7-JS, -PM$^8_{6,1}$, -RM$^7_{6,2,0}$ and -AIM$^5_{4,2;1E4}$.

Based on above results, the evolutions of L$_1$-error of schemes with $\Delta x$ are derived and shown in Fig. 20. The figure shows that WENO7-PRM$^4_{3,3}$ indicates relatively lager errors on the first two grids but achieve almost the same error on the third grids, while the rest schemes perform similarly. According to above discussions, the larger error of WENO7-PRM$^4_{3,3}$ on the first grids should owe to the deviation at the left foot of the fist peak (see Fig. 17(b)), while the error on the second grids owe to smeared left description of the first peak and additional oscillations between the first two structures (see Fig. 18(b)). It seems that WENO7-AIM$^5_{4,2;1E4}$ indicates a performance with overall less osccillations, however such perfomance does not guaranttee its robustness in subsequent problems by Euler equations.

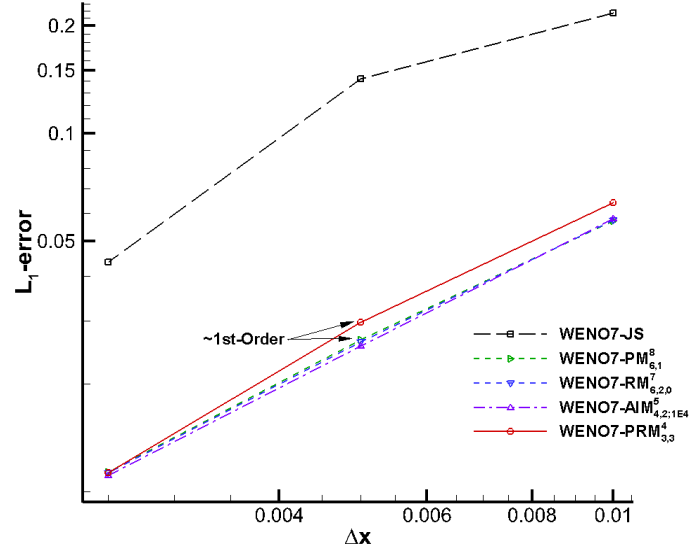

Fig. 20. $L_1$-errors with $\Delta x$ from computations of combination-waves advection at $T$=2000 by WENO7-PRM$^4_{3,3}$ with the comparisons by WENO7-JS, -PM$^8_{6,1}$,-RM$^7_{6,2,0}$ and -AIM$^5_{4,2;1E4}$

5.4 1-D problems by Euler equations

In this situation, the third-order TVD Runge-Kutta [2] is used for temporal discretization and Steger-Warming scheme is employed for flux splitting.

(1) Results of WENO3-PRM$^2_{1,1}$ and corresponding third-order comparatives

(a) Strong shock wave

Besides WENO3-PRM$^2_{1,1}$, the following mapped WENO3-JS schemes with $r_{c,g}$=2 are chosen for comparison, namely WENO3-M, -PPM$^3_{2,0}$ and -IM$^3_{2,0;0.1}$; in addition, two recently proposed WENO-Z-type scheme are selected also: WENO3-P+3[14] and WENO3-F3[15]. In two computations with $PR$=10$^3$ and 10$^6$, all schemes *except* WENO3-P+3 fullfill the computation. Because all results are similar in the first case, they are omitted for brevity, and the results of the second case are illustrated in Fig. 21. The figure tells that all schemes do not well simulate the platform of density after the shock, where WENO3-F3 seems to yield a result with relatively lower peak. Because WENO3-M shows a peak with higher height than that of WENO3-IM$^3_{2,0;0.1}$, we would rather regard this case as a test of numerical robustness than indication of resolution.

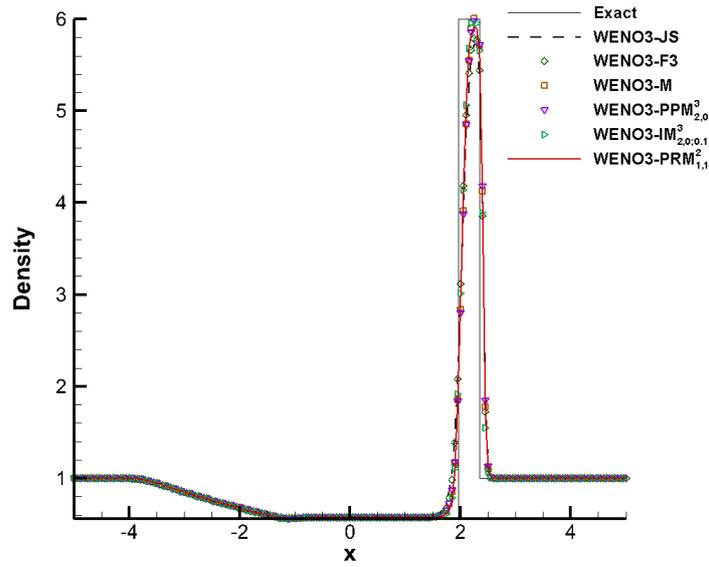

Fig. 21. Density distribution of strong shock wave with initial pressure ratio $PR=10^6$ at $t=0.01$ by WENO3-PRM$^2_{1,1}$ with the comparisons by WENO3-F3, -M, -PPM$^3_{2,0}$ and -IM$^3_{2,0;0.1}$.

(b) Blast wave

Besides WENO3-PRM$^2_{1,1}$, similar third-order comparatives are chosen as above with the density distributions shown in Fig. 22. It is found that all schemes can fulfill the computation except WENO3-P+3. According to our experience, the third-order schemes ususally perform robustly because of their small stencils, therefore the failure of WENO3-P+3 indicates even the third-order schemes should not only focus on resolution but also the numerical stability. Considering the performance of WENO3-PRM$^2_{1,1}$ later on Shu-Osher problem, the scheme indicates a well balance between the robustness and resolution.

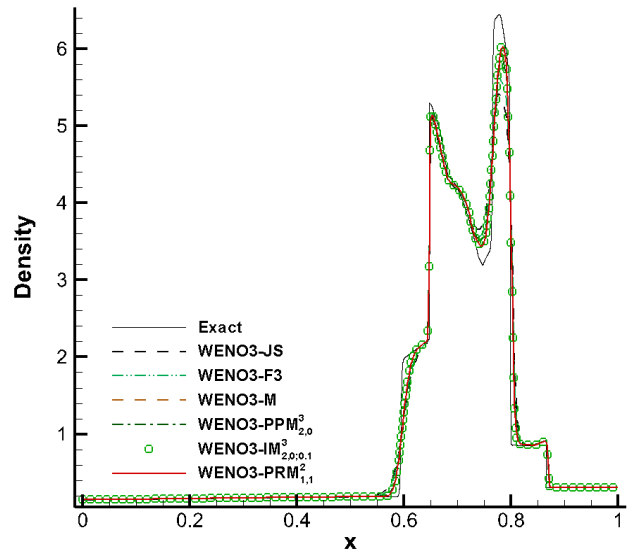

Fig. 22. Density distributions of blast waves at $t=0.038$ on 200 grids by WENO3-PRM$^2_{1,1}$ with the comparisons by WENO3-F3, -M, -PPM$^3_{2,0}$ and -IM$^3_{2,0;0.1}$.

(c) Shu-Osher problem

Besides WENO3-PRM$^2_{1,1}$, similar comparatives are adopted for testing as in previous problem. Other than grids with number usually from 400-600, the employment of 240 grids provides as a tough test on numerical resolution. Corresponding results are shown in Fig. 23. One can see that only the first three waves after the shock wave are distinguishable by WENO-PRM$^2_{1,1}$, -IM$^3_{2,0;0.1}$ and -P+3. Among the three schemes, WENO-PRM$^2_{1,1}$ and -IM$^3_{2,0;0.1}$ yield results with the best resolution. What is more, only WENO-PRM$^2_{1,1}$ tells the density fluctuation of the fourth wave structure and indicates its superior resolution thereby.

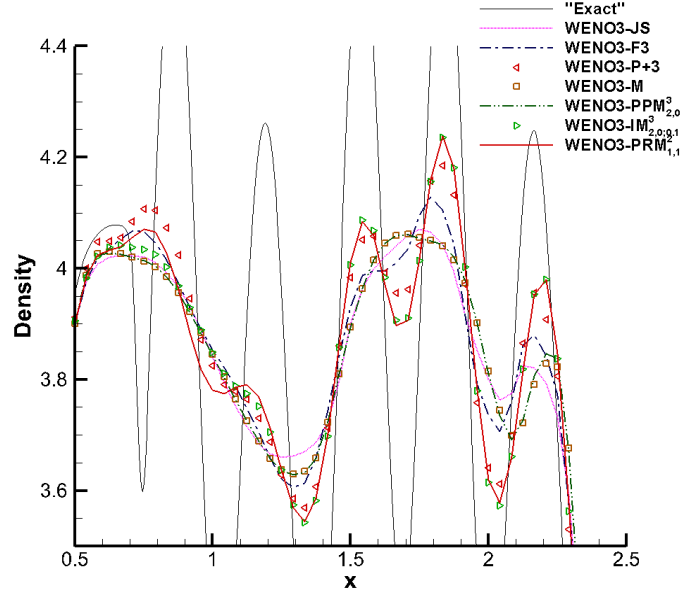

Fig. 23. Density distributions of Shu-Osher problem at $t$=1.8 on 240 grids by WENO3-PRM$^2_{1,1}$ with the comparisons by WENO3-JS, WENO3-P+3, -F3, -M, -PPM$^3_{2,0}$ and -IM$^3_{2,0;0.1}$.

(2) Results of WENO5-PRM$^3_{2,2}$ and corresponding fifth-order comparatives

(a) Strong shock wave

To compare with PRM$^3_{2,2}$, the following representative mappings in fifth-order scenario are chosen: $g_M$, PM$^8_{6,1}$, IM$^3_{2,0;0.1}$ and RM$^7_{6,2,0}$, where the latter three are intensively studied in Refs. 7, 8 and 10 as improvements for WENO5-M. Besides comparative mappings, some WENO-Z type schemes are chosen as well, namely WENO-Z with $q$=1 and 2 [4] and WENO-NIS [17]. All schemes can accomplish the computations with density distributions shown in Fig. 24. One can see that the theoretical platform after the shock is still not well resolved. Specifically, WENO5-Z at $q$=1 and WENO5-IM$^3_{2,0;0.1}$ yield the relatively larger overshoot over the exact solution while WENO5-PRM $^3_{2,2}$ inidcates a moderate distribution.

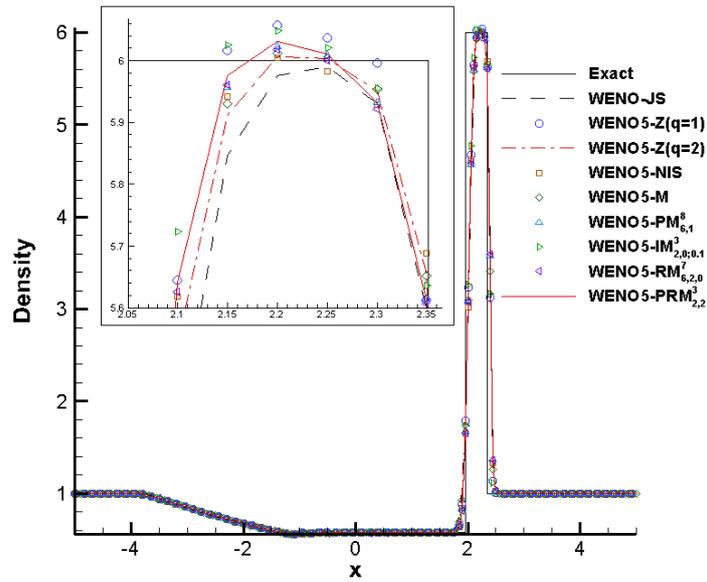

Fig. 24. Density distribution of strong shock wave with initial pressure ratio $PR=10^6$ at $t=0.01$ by WENO5-PRM$^3_{2,2}$ with the comparisons by WENO5-Z($q=1$ and 2), -NIS, -M, -PM$^8_{6,1}$, -IM$^3_{2,0;0.1}$ and -RM$^7_{6,2,0}$.

(b) Blast wave

To compare with WENO5-PRM$^3_{2,2}$, similar fifth-order comparatives are chosen as above. The density distributions are shown in Fig. 25 except that of WENO5-Z at $q=1$ and -IM$^3_{2,0;0.1}$. The two schemes, which have yielded the relatively largest overshoot in strong shock wave problem, fail to accomplish the computation because of the blow-up. According to Ref. 4, WENO5-Z at $q=1$ would yield higher resolution than at $q=1$, however, the potential less robustness exists as shown here.

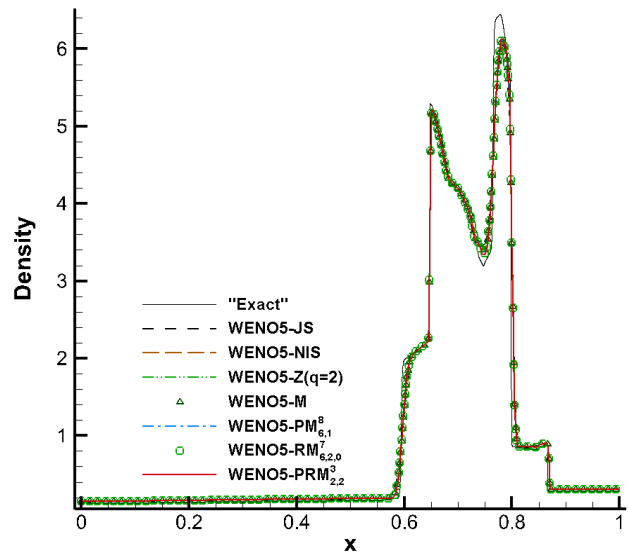

Fig. 25. Density distributions of blast waves at $t=0.038$ on 200 grids by WENO5-PRM$^3_{2,2}$ with the comparisons by WENO5-NIS, -Z($q=2$), -M, -PM$^8_{6,1}$ and -RM$^7_{6,2,0}$.

(c) Shu-Osher problem

To compare with WENO5-PRM$^3_{2,2}$, similar comparatives are chosen for testing as just now. For calrity, the results of comparatives are shown in two groups, i.e. results of mapped WENO5 in Fig. 26(a) and that of WENO-Z type schemes in Fig. 26(b). From Fig. 30(a), one can see that WENO5-PRM$^3_{2,2}$ and -IM$^3_{2,0;0.1}$ indicate the best resolution, where the latter even displays slightly better capability; WENO-RM$^7_{6,2,0}$ shows a resolution of the second class, while the rest schemes exhibit a relative inferior resolution. Recalling that WENO-IM$^3_{2,0;0.1}$ fails in blast wave computation, one can see that WENO5-PRM$^3_{2,2}$ demosntrate outstanding performance with well balance between the resolution and robustness among mapped WENO5 schemes.

Nest, the result of WENO5-PRM$^3_{2,2}$ is compared with that of WENO-Z type schemes in Fig. 26(b). The figure tells that WENO5-PRM$^3_{2,2}$ and WENO5-Z at $q=1$ achieves the best resolution and outperform the rest WENO-Z schemes. Similarly, recalling the failure of WENO5-Z at $q=1$ in blast wave computation, the comparative advantage of WENO5-PRM$^3_{2,2}$ is manifested once again.

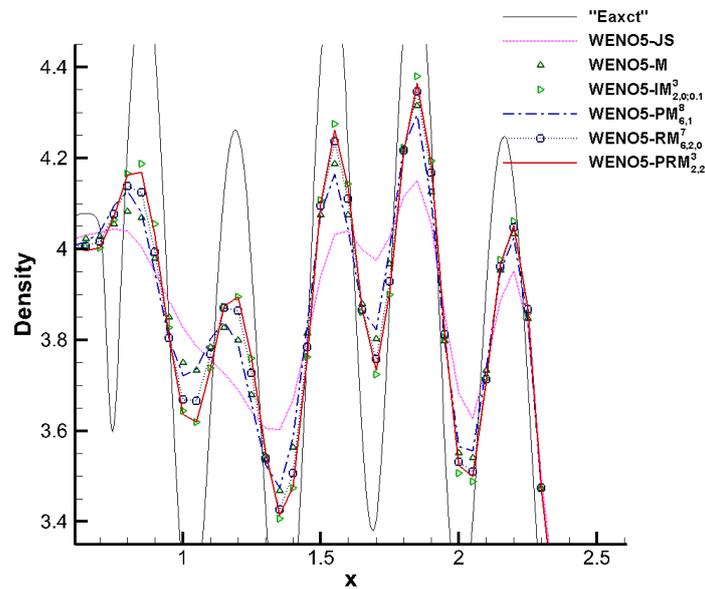

(a) Results compared by WENO5-JS, -M, -IM$^3_{2,0;0.1}$, -PM$^8_{6,1}$ and -RM$^7_{6,2,0}$

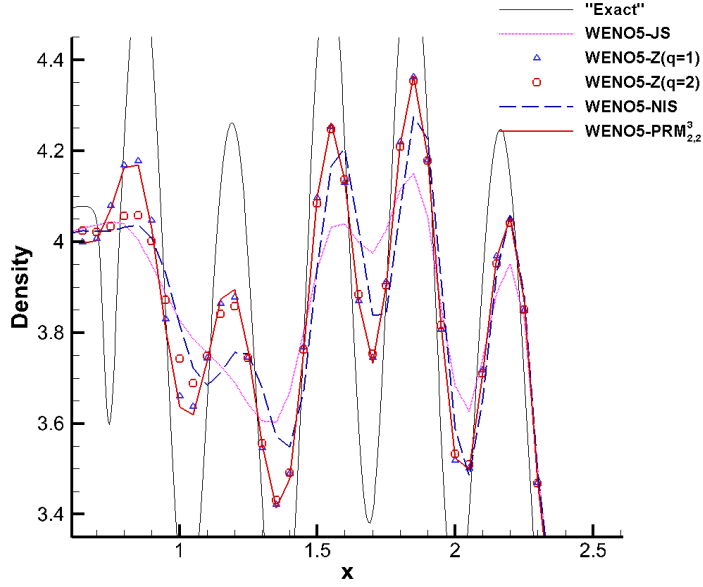

(b) Results compared by WENO5-JS, -NIS and -Z($q$=1, 2)

Fig. 26. Density distributions of Shu-Osher problem at $t$=1.8 on 200 grids by WENO5-PRM$^3_{2,2}$ with the comparisons by WENO5-JS, other mapped WENO5 and WENO5-Z type schemes.

(d) Titarev-Toro problem

To compare with WENO5-PRM$^3_{2,2}$, similar mapped schemes and WENO-Z type methods are employed. Because of the high frequency of the density fluctuation, only parts of the distributions are displayed where possible differences within results of schemes are visible, i.e. the region at $x \in$[1.33, 3.5] as shown in Fig. 27(a), whilst zoomed views are further displayed subsequently. The area just after the shock is first zoomed in Fig. 27(b), where the fluctuations undergo numerical dissipation shortly; afterwards, a further suffering of dissipation at $x \in$[1.36, 2.12] is zoomed in Fig.27(c). The figures tell that WENO5-PRM$^3_{2,2}$, -IM$^3_{2,0;0.1}$ and WENO5-Z with $q$=2 achieve the best resolution; WENO5-RM$^7_{6,2,0}$ and WENO5-Z with $q$=1 take the seond place; WENO5-M and -PM$^8_{6,1}$ take the third position where the former inidcates relatively better resolution. Note $r_{c,g}$ in PM$^8_{6,1}$ is 6, much larger than 2 in PRM$^3_{2,2}$ and IM$^3_{2,0;0.1}$, therefore the mappings with samller $r_{c,g}$ but large flatness can outperform the one with larger $r_{c,g}$ but small flatness. WENO5-NIS fails in the computation, which reminds once again the importance of robustness in developing high-order schemes.

Another point worthy of attention is that WENO5-Z with $q$=2 slightly outperforms the scheme at $q$=1, which is oppsite to the resolution relationship shown in Shu-Osher problem (see Fig. 26(b)). In Ref. 4, WENO5-Z at $q$=1 was thought to be one order lower than the scheme at $q$=2; however, the former would allocate more weights to unsmooth stencils and higher resolution was expected usually. The distributions in Fig. 27 indicate the advantage of accuracy order would somtimes exceed effects of nonliear techniques. Unsurprisingly, WENO5-JS indicates a poor resolution and almost fails to tell downstream fluctuations,which is the same as that in Refs. 10-11.

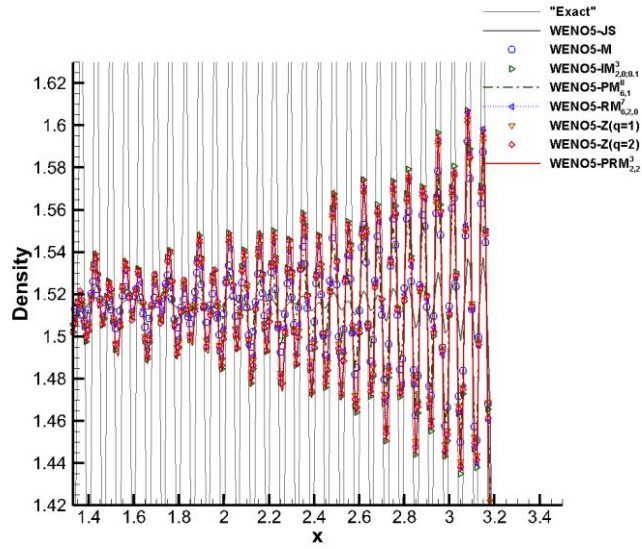

(a) Global view where $x \in [1.33, 3.5]$

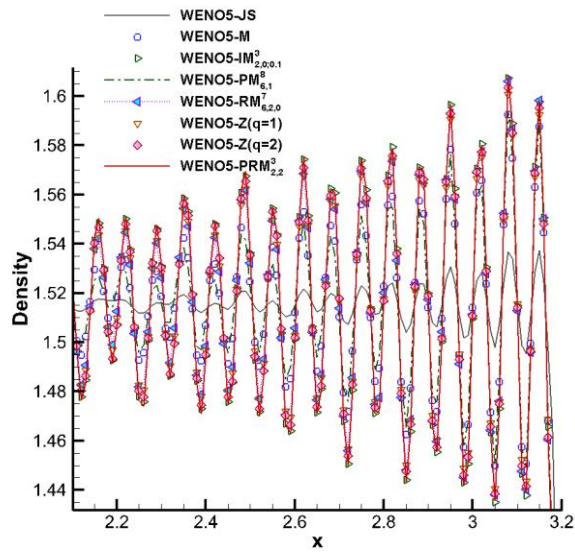

(b) Zoomed view 1 where $x \in [2.1, 3.2]$

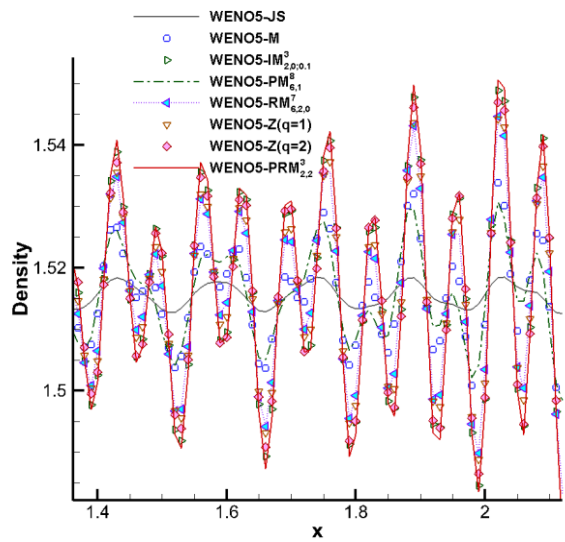

(c) Zoomed view 2 where $x \in [1.36, 2.12]$

Fig. 27. Local and zoomed view of density distribution of Titarev-Toro problem at $t=5$ on 1000 grids by WENO5-PRM$^3_{2,2}$ with the comparisons by WENO5-JS, -M, -IM$^3_{2,0;0.1}$, -PM$^8_{6,1}$, -RM$^7_{6,2,0}$, and WENO5-Z with $q=1$ and 2

(3) Results of WENO7-PRM$^4_{3,3}$ and corresponding seventh-order comparatives

(a) Strong shock wave

As in previous scalar problems, the similar mappings are chosen to compare with PRM$^4_{3,3}$: PM$^8_{6,1}$, RM$^7_{6,2,0}$ and AIM$^5_{4,2;1E4}$-M. The use of AIM$^5_{4,2;1E4}$-M other than AIM$^5_{4,2;1E4}$ is to follow the suggestion of Ref. 11 that the former will dehave more robustly. And as usual, the result of canonical WENO7-JS is provided. It is surprising to note that WENO7-AIM$^5_{4,2;1E4}$-M fails to accomplish the computation at both pressure ratios because of the blow up, which implies the less robustness despite of its well performance in previous scalar computations. The other schemes fufill the test, and the results at $PR=10^6$ are shown in Fig. 28. The figure shows that the schemes yield similar results in which the density platform is still not well resolved.

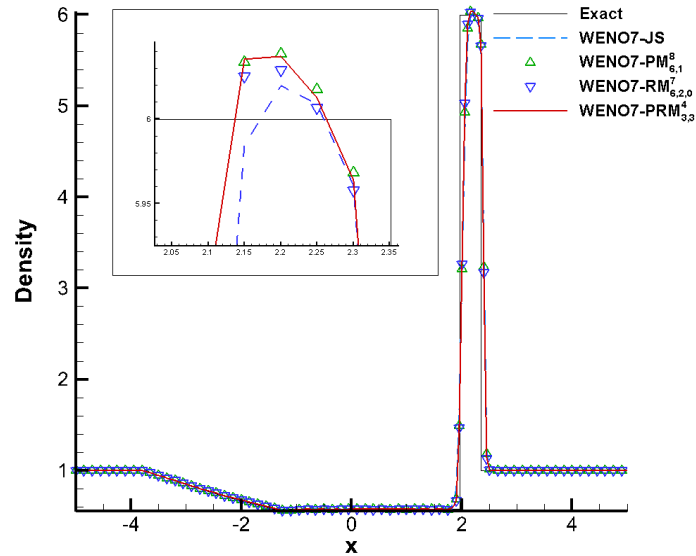

Fig. 28. Density distribution of strong shock wave with initial pressure ratio $PR=10^6$ at $t=0.01$ by WENO7-PRM$^4_{3,3}$ with the comparisons by WENO7-JS, -PM$^8_{6,1}$ and -RM$^7_{6,2,0}$.

(b) Blast wave

Besides WENO7-PRM$^4_{3,3}$, similar comparative schemes are chosen as above. Again, WENO7-AIM$^5_{4,2;1E4}$-M fails in the computation because of the blow up. The other schemes yield similar results which are shown in Fig. 29 with details in zoomed view.

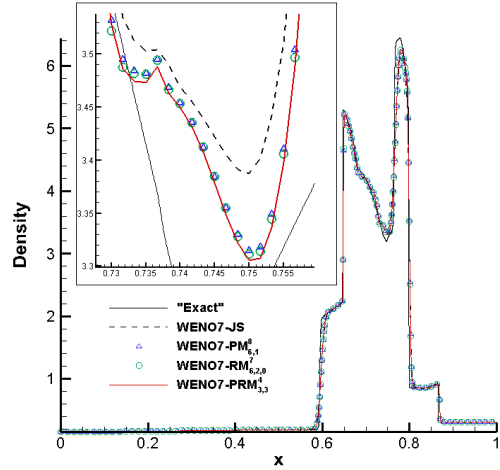

Fig. 29. Density distributions of blast waves at $t$=0.038 on 200 grids by WENO7-PRM$^4_{3,3}$ with the comparisons by WENO7-JS, -PM$^8_{6,1}$ and -RM$^7_{6,2,0}$.

(c) Shu-Osher problem

To compare with WENO7-PRM$^4_{3,3}$, the similar comparatives are adopted for testing as before, and corresponding results are shown Fig. 30. The figure tells that WENO7-PRM$^4_{3,3}$, RM$^7_{6,2,0}$ and AIM$^5_{4,2;1E4}$-M yield results with almost the same resoltuion, whereas the result of WENO7-PM$^8_{6,1}$ appears slightly dissipated.

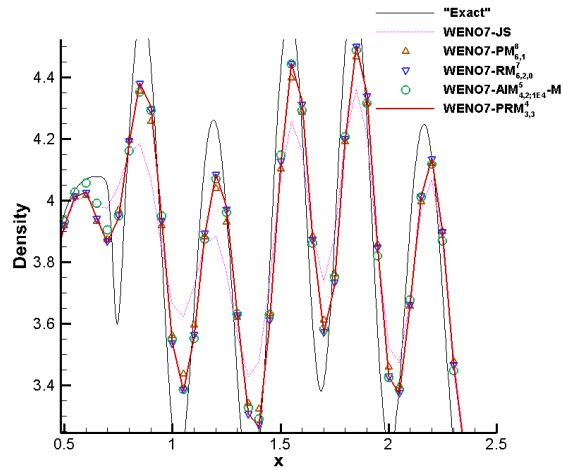

Fig. 30. Density distributions of Shu-Osher problem at $t$=1.8 on 200 grids by WENO7-PRM$^4_{3,3}$ with the comparisons by WENO7-JS, -PM$^8_{6,1}$, -RM$^7_{6,2,0}$ and -AIM$^5_{4,2;1E4}$-M.

(d) Titarev-Toro problem

Besides WENO7-PRM$^4_{3,3}$, the similar comparatives are chosen for testing as above, and corresponding results in local region are shown in Fig. 31. First, the overview of density fluctuation after the shock is shown in Fig. 31(a). To illustrate the details, the distributions in Fig. 31(a) is separately displayed in three zoomed views shown in Figs. 31(b)-(d). One can see that WENO7-PRM$^4_{3,3}$, performs almost the same as WENO7-RM$^7_{6,2,0}$, while WENO7-AIM$^5_{4,2;1E4}$-M inidcates slightly higher resultion. By contrast, WENO7-PM$^8_{6,1}$ yields relatively less-resolved description of fluctuations, although the scheme has resolved the structures as well. As the comparison, WENO7-JS indicates poor resolution on the structures and even fails to resolve fluctuations in the downstrean region as shown in Fig. 31(c)-(d).

Although WENO7-AIM$^5_{4,2;1E4}$-M inidicates slight superiority in this case, its failure in strong shock wave and blast wave problems reminds sufficient attention should be paid in the pursue of high resolution.

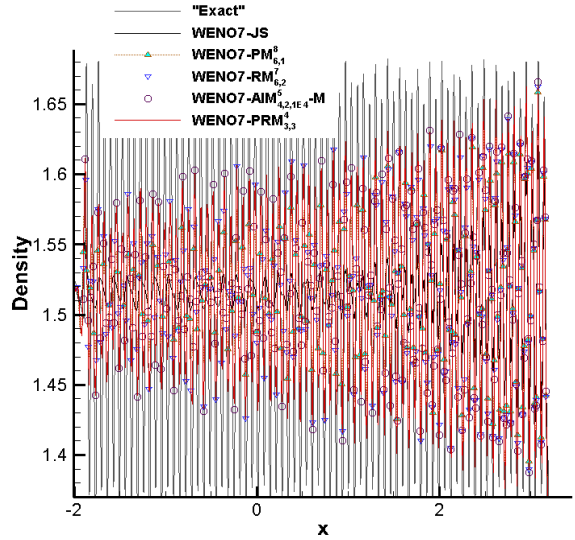

(a) Local view where $x \in [1.33, 3.5]$

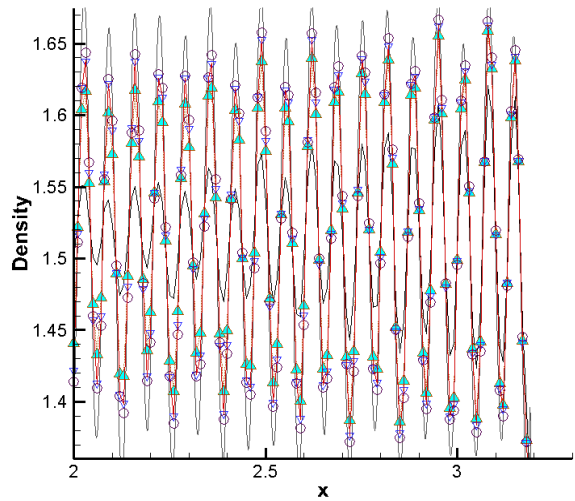

(b) Zoomed view 1 where $x \in [2, 3.3]$

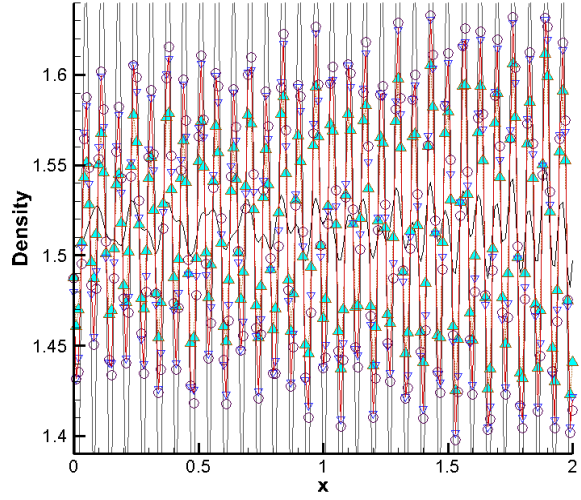

(c) Zoomed view 2 where $x \in [0, 2]$

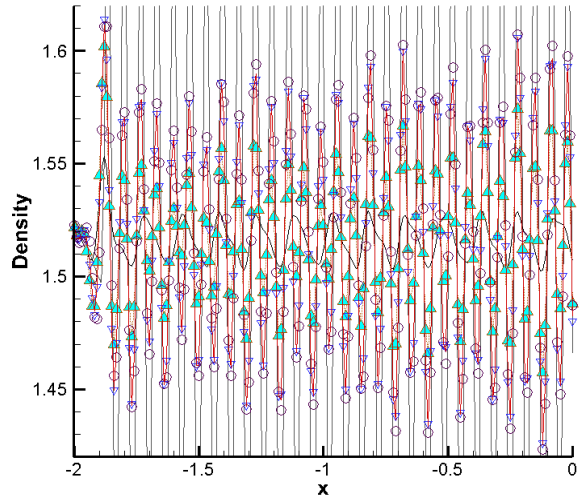

(d) Zoomed view 2 where $x \in [-2, 0]$

Fig. 31. Local and zoomed view of density distribution of Titarev-Toro problem at $t=5$ at 1000 grids by by WENO7-PRM$^4_{3,3}$ with the comparisons by WENO7-JS, -PM$^8_{6,1}$, -RM$^7_{6,2,0}$ and AIM$^5_{4,2;1E4}$-M.

**6 Conclusions**

Comprehensive and intensive investigations are carried out regarding piecewise rational mapping, and corresponding new PRM method is proposed with sufficient regulation capability. Through theoretical analysis and numerical tests, the following conclusions are drawn:

(1) So-called $C_{n,m}$ condition is summarized to develop rational mapping, which incarnates the favorable properties of an ideal mapping. In $C_{n,m}$, the endpoint convergence pattern as WENO, which was first proposed by Ref. 9, is utilized in this paper. On the one hand numerical example indicates the pattern as WENO performs similarly as that as ENO [6], on the other hand the former has no restrictions on the choice of $n$ and $m$ of $C_{n,m}$ other than the latter does.

(2) Numerical examples indicate that the piecewise implementation of mapping does not necessarily entail or favor numerical instability in high-order mapped WENO-JS (e.g. mapped WENO9-JS), which is contrary to the conclusion in Ref. 10.

(3) By upgrading the smoothness indicators of WENO3-JS to that of WENO5-JS as Ref. 13, WENO3-JS can preserve its optimal order through mapping in the occurrence of first-order critical points, which is thought to be unavailable in Refs. 3, 7-8 and 10-11. Careful analysis is provided and numerical validations are shown in this regard.

(4) New piecewise rational mapping $PRM^{n+1}_{n,m}$ is proposed. The mapping on the one hand can explicitly make the desired flatness and endpoint convergence properties achievable simultaneously, on the other hand exhibits powerful capability of regulation. As an example for the latter, one particular $PRM^3_{2,2}$ can resemble existing mappings, namely $PM^8_{6,1}$ [7] and $RM^7_{6,2,0}$ [10] by specific choices of parameters. Two theorems are provided with proofs to guarantee $PRM^{n+1}_{n,m}$ satisfy $C_{n,m}$ and to define the valid range of parameters of mapping. As a byproduct, the general form of $PPM^{n+m+1}_{n,m}$ is proposed and one theorem is provided with proof to guarantee $C_{n,m}$ satisfied.

Three concrete ones, i.e. $PRM^2_{1,1}$, $PRM^3_{2,2}$ and $PRM^4_{3,3}$ are determined for WENO3,5,7-JS, and corresponding WENO3-$PRM^2_{1,1}$, WENO5-$PRM^3_{2,2}$, WENO7-$PRM^4_{3,3}$ are acquired.

(5) Regarding WENO3-$PRM^2_{1,1}$, although $PRM^2_{1,1}$ does not have a larger Flatness-I than $g_M$ and $IM^3_{2,0;0.1}$ in the neighborhood of linear weights, it has the largest Flatness-II among comparatives as well as the best performance regarding *PEC*. Numerically, WENO3-$PRM^2_{1,1}$ not only can preserve the third-order accuracy in the occurrence of first-order critical points but also shows the quickest rate of order convergence among comparatives. In computations regarding robustness and stability, WENO3-$PRM^2_{1,1}$ shows one of the best performances whereas WENO3-P+3 fails in tests of strong shock and blast wave; besides, WENO3-$PRM^2_{1,1}$ shows the best resolution in Shu-Osher problem on 240 grids which even surpasses that of WENO-$IM^3_{2,0;0.1}$. Hence WENO3-$PRM^2_{1,1}$ indicates the comprehensive advantages in accuracy, robustness and resolution.

(6) Regarding WENO5-$PRM^3_{2,2}$, $PRM^3_{2,2}$ shows excellent Flatness-II and *PEC* in balance among comparatives such as $g_M$, $IM^3_{2,0;0.1}$, $PM^8_{6,1}$ and $RM^7_{6,2,0}$; moreover, the drawbacks such as excessive flatness and poor performance of *PEC* of the latter three are overcome. Numerically, WENO5-$PRM^3_{2,2}$ shows one of the quickest rates of order convergence, while preserving the fifth-order accuracy at first-order critical points. In regard to stability of long time computation, WENO5-$PRM^3_{2,2}$ indicates the similar performance as that of $PM^8_{6,1}$ and $RM^7_{6,2,0}$. Although WENO5-$IM^3_{2,0;0.1}$ shows quite slight improvement, it fails in blast wave problem, which indicates aforementioned stability focused by Ref. 7-8 and 10 does not equal to robustness. A test $R^3_{2,2}$ which has similar Flatness-II as that of $PRM^3_{2,2}$ produces oscillations in long time computation, which verifies the necessity of extra implementation in $PRM^3_{2,2}$ to guarantee numerical stability by sufficient regulation. In computations regarding robustness and stability, WENO5-$PRM^3_{2,2}$ shows one of the best performances in presence of failures of WENO5-$IM^3_{2,0;0.1}$ and WENO5-Z at $q=1$ in blast wave. Similarly, WENO5-$PRM^3_{2,2}$ shows one of the best resolution in Shu-Osher and Titarev-Toro problems, whereas the tiny improvement by WENO5-$IM^3_{2,0;0.1}$ reminds the risk of impaired robustness if only focusing resolution.

(7) Regarding WENO7-$PRM^4_{3,3}$, $PRM^4_{3,3}$ shows well-performed Flatness-II and *PEC* in balance among comparatives such as $PM^8_{6,1}$ and $RM^7_{6,2,0}$. Numerically, WENO7-$PRM^4_{3,3}$ shows similar quick rate of order convergence with the preservation of seventh-order accuracy at first- and second-order critical points in spite of relatively larger errors on first two grids. In long time computation, WENO7-PRM$^4_{3,3}$ indicates similar performance as that of WENO7-PM$^8_{6,1}$ and -RM$^7_{6,2,0}$, whereas WENO7-AIM$^5_{4,2;1E4}$ inidcates tiny improvements. Similar situations occur in the computations of Shu-Osher and Titarev-Toro problems; however, WENO7-AIM$^5_{4,2;1E4}$(-M) indicates an obvious inferior robustness because of its failure in strong shock and blast wave problems. In summary, according to numerical tests conducted, WENO7-PRM$^4_{3,3}$ indicates a similar performance as that of WENO7-RM$^7_{6,2,0}$, while the profile of PRM$^4_{3,3}$ shows more advantage on stability than that of RM$^7_{6,2,0}$ due to performances near endpoints. Additionally, AIM$^5_{4,2;1E4}$ is found to fail to work with WENO5-JS.

**Acknowledgements**

This study is sponsored by the project of National Numerical Wind-tunnel of China under the grant number NNW2019ZT4-B12.

**Appendix I**

For reference, the coefficients of candidate schemes, linear weights, and coefficients of smoothness indicators of WENO-JS which correspond to Eqns. (3)-(5) are tabulated in Tables 7-11.

Table 8 Coefficients $a^r_{kl}$ of candidate schemes $q^r_k$ of WENO-JS with $r$=2~5 [2]

| r | k | $a^r_{k0}$ | $a^r_{k1}$ | $a^r_{k2}$ | $a^r_{k3}$ | $a^r_{k4}$ |
|---|---|---|---|---|---|---|
| 2 | 0 | -1/2 | 3/2 | - | - | - |
|   | 1 | 1/2 | 1/2 | - | - | - |
| 3 | 0 | 2/6 | -7/6 | 11/6 | - | - |
|   | 1 | -1/6 | 5/6 | 2/6 | - | - |
|   | 2 | 2/6 | 5/6 | -1/6 | - | - |
| 4 | 0 | -3/12 | 13/12 | -23/12 | 25/12 | - |
|   | 1 | 1/12 | -5/12 | 13/12 | 3/12 | - |
|   | 2 | -1/12 | 7/12 | 7/12 | -1/12 | - |
|   | 3 | 3/12 | 13/12 | -5/12 | 1/12 | - |
| 5 | 0 | 12/60 | -63/60 | 137/60 | -163/60 | 137/60 |
|   | 1 | -3/60 | 17/60 | -43/60 | 77/60 | 12/60 |
|   | 2 | 2/60 | -13/60 | 47/60 | 27/60 | -3/60 |
|   | 3 | -3/60 | 27/60 | 47/60 | -13/60 | 2/60 |
|   | 4 | 12/60 | 77/60 | -43/60 | 17/60 | -3/60 |

Table 9 Linear weights $d^r_k$ of WENO-JS with $r$=2~5 [2]

| r | $d_0$ | $d_1$ | $d_2$ | $d_3$ | $d_4$ |
|---|---|---|---|---|---|
| 2 | 1/3 | 2/3 | - | - | - |
| 3 | 1/10 | 6/10 | 3/10 | - | - |
| 4 | 1/35 | 12/35 | 18/35 | 4/35 | - |
| 5 | 1/126 | 20/126 | 60/126 | 40/126 | 5/126 |

Table 10 Coefficients $b^r_{kml}$ of smoothness indicators of WENO-JS in Eq. (5) with $r$=2~5 [2, 10-12]

| r | k | m | $b^r_{km0}$ | $b^r_{km1}$ | $b^r_{km2}$ | $b^r_{km3}$ | $b^r_{km4}$ | k | m | $b^r_{km0}$ | $b^r_{km1}$ | $b^r_{km2}$ | $b^r_{km3}$ | $b^r_{km4}$ |
|---|---|---|---|---|---|---|---|---|---|---|---|---|---|---|
| 2 | 0 | 0 | -1 | 1 | - | - | - | | | | | | | |
| | 1 | 0 | -1 | 1 | - | - | - | | | | | | | |
| 3 | 0 | 0 | 1 | -4 | 3 | - | - | 2 | 0 | 3 | -4 | 1 | - | - |
| | | 1 | 1 | -2 | 1 | - | - | | 1 | 1 | -2 | 1 | - | - |
| | 1 | 0 | -1 | 0 | 1 | - | - | | | | | | | |
| | | 1 | 1 | -2 | 1 | - | - | | | | | | | |
| 4 | 0 | 0 | -2 | 9 | -18 | 11 | - | 2 | 0 | -2 | -3 | 6 | -1 | - |
| | | 1 | -1 | 4 | -5 | 2 | - | | 1 | 1 | -2 | 1 | 0 | - |
| | | 2 | -1 | 3 | -3 | 1 | - | | 2 | -1 | 3 | -3 | 1 | - |
| | 1 | 0 | 1 | -6 | 3 | 2 | - | 3 | 0 | -11 | 18 | -9 | 2 | - |
| | | 1 | 0 | 1 | -2 | 1 | - | | 1 | 2 | -5 | 4 | -1 | - |
| | | 2 | -1 | 3 | -3 | 1 | - | | 2 | -1 | 3 | -3 | 1 | - |
| 5 | 0 | 0 | 3 | -16 | 36 | -48 | 25 | 3 | 0 | 3 | 10 | -18 | 6 | -1 |
| | | 1 | 119 | -606 | 1234 | -1126 | 379 | | 1 | 119 | -216 | 64 | 44 | -11 |
| | | 2 | 3 | -14 | 24 | -18 | 5 | | 2 | 3 | -10 | 12 | -6 | 1 |
| | | 3 | 1 | -4 | 6 | -4 | 1 | | 3 | 1 | -4 | 6 | -4 | 1 |
| | 1 | 0 | 1 | -6 | 18 | -10 | -3 | 4 | 0 | 25 | -48 | 36 | -16 | 3 |
| | | 1 | 11 | -44 | -64 | 216 | -119 | | 1 | 379 | -1126 | 1234 | -606 | 119 |
| | | 2 | 1 | -6 | 12 | -10 | 3 | | 2 | 5 | -18 | 24 | -14 | 3 |
| | | 3 | 1 | -4 | 6 | -4 | 1 | | 3 | 1 | -4 | 6 | -4 | 1 |
| | 2 | 0 | 1 | -8 | 0 | 9 | -1 | | | | | | | |
| | | 1 | 11 | -174 | 326 | -174 | 11 | | | | | | | |
| | | 2 | 1 | -2 | 0 | 2 | -1 | | | | | | | |
| | | 3 | 1 | -4 | 6 | -4 | 1 | | | | | | | |

Table 11 Coefficients $c^r_m$ of smoothness indicators of WENO-JS in Eq. (5) with $r$=2~5 [2, 10-12]

| r | $c^r_0$ | $c^r_1$ | $c^r_2$ | $c^r_3$ |
|---|---|---|---|---|
| 2 | 1 | - | - | - |
| 3 | 1/4 | 13/12 | - | - |
| 4 | 1/36 | 13/12 | 781/720 | - |
| 5 | 1/144 | 13/202800 | 781/2880 | 1421461/1310400 |

Considering the occurrence of critical point with the order $n_{cp}$, suppose $r_c$-WENO-JS is corresponding order of WENO-JS, $r_c$-WENO-M is that of WENO-M, and $r_{c,g}$ is the maximum order of $g^{(i)}(\omega_k)=0$ for order recovery of WENO-JS (see Eq. (8)). According to Refs. 3, the above orders can be tabulated in Table 12 for WENO-JS from $r$=2 to 5.

Table 12 Relationship of $r_c$-WENO-JS, $r_c$-WENO-M and $r_{c,g}$ with $r$ and $n_{cp}$

| | $n_{cp}$ | $r_c$-WENO-JS | $r_c$-WENO-M | $r_{c,g}$ |
|---|---|---|---|---|
| $r$=2 | 0 | 3 | - | - |
| | 1 | 1 | N/A | N/A |
| $r$=3 | 0 | 5 | - | - |
| | 1 | 3 | 5 | ≥2 |
| | 2 | 2 | 2 | N/A |

| r=4 | 0 | 7 | - | - |
|---|---|---|---|---|
| | 1 | 5 | 7 | ≥2 |
| | 2 | 4 | 6 | ≥3 |
| | 3 | 3 | 3 | N/A |
| r=5 | 0 | 9 | 9 | - |
| | 1 | 7 | 9 | ≥2 |
| | 2 | 6 | 9 | ≥2 |
| | 3 | 5 | 7 | ≥4 |
| | 4 | 4 | 4 | N/A |

It is worthwhile to mention that when the upgradation of smoothness indicator described in Section 3.1 is employed for WENO3, its optimal order could be recovered by $r_{c,g} \geq 1$.

**Appendix II**

Given the mapping $R_{n,m}^{L,n+1}$ in $[d_k, 1]$ in Eq. (11), the symmetric $R_{n,m}^{L,n+1}$ can be derived by $R_{n,m}^{L,n+1} = \frac{d_k}{1-d_k}\left[1 - R_{n,m}^{R,n+1}(1 - \frac{1-d_k}{d_k}\omega)\right]$ in the form as Eq. (12). Corresponding coefficients $c_{n,m,i}^{L}$ except $c_{n,m,4}^{L}$ are tabulated in Table 13 for $m, n \leq 5$. Almost all $c_{n,m,4}^{L}$ equals zero except $c_{4,2,4}^{L} = -d_k\left(\frac{b}{1-d_k} - 4\right)$.

Table 13 Coefficients $c_i^{L,n}$ except $c_4^{L,n}$ of $(n, i, m)$ in mapping $R_{n,m}^{L,n+1}$

| | | m=0 | m=1 | m=2 | m=3 | m=4 |
|---|---|---|---|---|---|---|
| n=1 | i=1 | 0 | 1 | - | - | - |
| | i=2 | b | $\frac{-b}{1-b \cdot d_k}$ | - | - | - |
| | i=3 | -1 | 0 | - | - | - |
| n=2 | i=1 | 0 | $-2d_k$ | 1 | - | - |
| | i=2 | $\frac{-b \cdot d_k}{1-d_k}$ | b | $\frac{b \cdot d_k - b}{d_k(1-b \cdot d_k)}$ | - | - |
| | i=3 | 1 | -1 | 0 | - | - |
| n=3 | i=1 | 0 | $3d_k^2$ | 1 | 1 | - |
| | i=2 | $\frac{b \cdot d_k^2}{(1-d_k)^2}$ | $\frac{-b \cdot d_k}{1-d_k}$ | b-1 | $\frac{b}{1-b \cdot d_k}$ | - |
| | i=3 | -1 | 2 | 0 | 0 | - |
| n=4 | i=1 | 0 | $-4d_k^3$ | 1 | 1 | 1 |
| | i=2 | $\frac{-b \cdot d_k^3}{(1-d_k)^3}$ | $\frac{b \cdot d_k^2}{(1-d_k)^2}$ | -1 | b-1 | $\frac{-b(1-d_k)}{d_k(1-b \cdot d_k)}$ |
| | i=3 | 1 | -3 | 0 | 0 | 0 |

**Appendix III**

In the following, proofs are provided for proposed theorems.

**Theorem 1**: Considering a function as $f(\omega) = d_k + \frac{1}{(1-d_k)^{n+m}}(\omega - d_k)^{n+1} \sum_{i=0}^{m} a_i^m (1-\omega)^{m-i}(1-d_k)^i$ in $[d_k, 1]$ where $a_i^m = \frac{\prod_{j=0}^{m-1-i}(n+j)}{(m-i)!}$ for $i < m$ and $a_m^m = 1$, then $f(\omega)$ satisfies $C_{n,m}$.

**Proof** It is trivial that $f(d_k) = d_k$ and $f(1) = 1$. Let $A = \frac{1}{(1-d_k)^{n+m}}$, $u = (\omega - d_k)^{n+1}$ and

$v = \sum_{i=0}^{m} a_i^m (1-\omega)^{m-i}(1-d_k)^i$, then $f(\omega) = d_k + Auv$. Using Leibnitz law,

$f(\omega)^{(k)} = A \sum_{i=0}^{k} C_k^i u^{(i)} v^{(k-i)}$. Subsequently, $f(\omega)^{(k)}$ at two endpoints will be discussed respectively.

(1) $f(d_k)^{(k)}$. Because $u^{(k)}(\omega) = \left[\prod_{i=0}^{k-1}(n+1-i)\right] \cdot (\omega - d_k)^{n-k+1}$ when $1 \leq k \leq n$, then

$f(d_k)^{(k)} = 0$ for the prescribed $k$; moreover, considering $u^{(n+1)}(d_k) = \prod_{i=0}^{k-1}(n+1-i) \neq 0$ and

$v(d_k) = \sum_{i=0}^{m} a_i^m (1-d_k)^{m-i}(1-d_k)^i \neq 0$, $f(d_k)^{(n+1)} = \left(Au^{(n+1)}v\right)_{\omega=d_k} \neq 0 \ 1 \leq k \leq m$.

(2) $f(1)^{(k)}$. Because $v^{(k)}(\omega) = (-1)^k \sum_{i=0}^{m-k} \frac{\prod_{j=0}^{m-1-i}(n+j)}{(m-k-i)!}(1-\omega)^{m-k-i}(1-d_k)^i$ when $1 \leq k \leq m$, then

$v^{(k)}(1) = (-1)^k \prod_{j=0}^{k-1}(n+j)(1-d_k)^{m-k}$ for the prescribed $k$. Consequently, one can find that

$f'(1) = A(u'v + uv')|_{\omega=1} = 1$. When $1 < k \leq m$,

$f(1)^{(k)} = A \sum_{p=0}^{k} C_k^p u^{(p)} v^{(k-p)}\Big|_{\omega=1} = A\left(uv^{(k)} + \sum_{p=1}^{k-1} C_k^p u^{(p)} v^{(k-p)} + u^{(k)}v\right)\Big|_{\omega=1}$

$= \frac{1}{(1-d_k)^{n+m}} \begin{bmatrix} (1-d_k)^{n+1}(-1)^k \prod_{j=0}^{k-1}(n+j)(1-d_k)^{m-k} + \\ \sum_{p=1}^{k-1} C_k^p \prod_{i=0}^{p-1}(n+1-i)(1-d_k)^{n-p+1}(-1)^{k-p} \prod_{j=0}^{k-p-1}(n+j)(1-d_k)^{m-k+p} + \\ \prod_{i=0}^{k-1}(n+1-i)(1-d_k)^{n-k+1}(1-d_k)^m \end{bmatrix}$

$= \frac{1}{(1-d_k)^{k-1}}\left[(-1)^k \prod_{j=0}^{k-1}(n+j) + \sum_{p=1}^{k-1} C_k^p \prod_{i=0}^{p-1}(n+1-i)(-1)^{k-p} \prod_{j=0}^{k-p-1}(n+j) + \prod_{i=0}^{k-1}(n+1-i)\right]$

Through symbolic operation, the term within the bracket turns out to be zero for given $k \geq 2$, which stands even at $k=m+1$, and therefore $f(1)^{(k)} = 0$. When $k = m+1$, it is trivial $v^{(m+1)}(\omega) = 0$, and

$$f(1)^{(m+1)} = A\sum_{p=0}^{m+1} C_k^p u^{(p)} v^{(m+1-p)}\Big|_{\omega=1}$$

$$= \frac{1}{(1-d_k)^{n+m}}\left[\begin{array}{l}\sum_{p=1}^{m} C_{m+1}^p \prod_{i=0}^{p-1}(n+1-i)(1-d_k)^{n-p+1}(-1)^{m+1-p}\prod_{j=0}^{m-p}(n+j)(1-d_k)^{p-1} + \\ \prod_{i=0}^{m}(n+1-i)(1-d_k)^{n-m}(1-d_k)^m\end{array}\right].$$

$$= \frac{1}{(1-d_k)^m}[\sum_{p=1}^{m} C_{m+1}^p \prod_{i=0}^{p-1}(n+1-i)(-1)^{m+1-p}\prod_{j=0}^{m-p}(n+j) + \prod_{i=0}^{m}(n+1-i)]$$

Comparing the last term in above $f(1)^{(m+1)}$ with the last term in $f(1)^{(k)}$ at $k=m+1$, one can find that $f(1)^{(m+1)} = \frac{1}{(1-d_k)^m}(-1)^m \prod_{j=0}^{m}(n+j) \neq 0$.

In short, $f(\omega)$ satisfies $C_{n,m}$ condition.

**Theorem 2**: Consider a mapping as $R_{n,m}^{n+1} = d_k + \frac{(\omega-d_k)^{n+1}}{(\omega-d_k)^n + c(1-\omega)^{m+1}}$ in [$d_k$, 1] with $0 < d_k < 1$ and $n, m \geq 1$. $c > 0$ is the sufficient and necessary condition for $R_{n,m}^{n+1}$ to be void of singularity.

**Proof** Other than targeting at the direct problem, we first study the solution of $c$ through which the denominator in $R_{n,m}^{n+1}$ has zero value. Consider $(\omega-d_k)^n + c(1-\omega)^{m+1} = 0$ at $\omega \in [d_k, 1]$, which actually defines a function $c(\omega)$ as $c = -\frac{(\omega-d_k)^n}{(1-\omega)^{m+1}}$. Because $\frac{dc(\omega)}{d\omega} = -\frac{(\omega-d_k)^{n-1}}{(1-\omega)^{m+2}}[n(1-\omega)+(m+1)(\omega-d_k)]$, then $\frac{dc(\omega)}{d\omega} > 0$ at $\omega \in [d_k, 1]$ for $0 < d_k < 1$ and $n, m \geq 1$. Hence the function implied in the solution is monotone and there is one-to-one correspondence between $c$ and $\omega$ at $\omega \in [d_k, 1]$. Considering $c(d_k) = 0$ and $c(1) = -\infty$, the range of $c$ turns out to be $(-\infty, 0]$. In other words, for any $c \in (-\infty, 0]$, there is a $\omega \in [d_k, 1]$ such that $(\omega-d_k)^n + c(1-\omega)^{m+1} = 0$. Hence the range of c such that $(\omega-d_k)^n + c(1-\omega)^{m+1} \neq 0$ at $\omega \in [d_k, 1]$ is $c > 0$.

**Theorem 3:** Consider a function as $f(\omega) = d_k + \frac{(\omega-d_k)^{n+1}}{(\omega-d_k)^n + c_2(\omega-d_k)^{n_1}(1-\omega)^{m_1} + c_1(1-\omega)^{m+1}}$ in [$d_k$, 1] where $n \geq 1$, and $c_1$, $c_2$ have the same sign. If $m \geq 0$ and $m_1 \geq 1$, then $f(\omega)$ satisfies $C_{n,\min(m,m_1-1)}$.

Prior to the proof of the theorem, the following lemma is proposed and proofed.

**Lemma** Supposing a function $f(\omega)$ which is defined as $f(\omega) = (\omega - a)^N \varphi(\omega)$ with $\varphi(a) \neq 0$ and $N \geq 1$, then $f(\omega)$ satisfies $f^{(i)}(a) = \begin{cases} 0, & 0 \leq i < N \\ N!\varphi(a), & i = N \end{cases}$.

**Proof** It trivial that $f^{(0)}(a) = 0$. Using Leibnitz law, $f^{(i)}(\omega) = \sum_{k=0}^{i} C_i^k \left[(\omega-a)^N\right]^{(k)} \varphi^{(i-k)}(\omega)$ with $0 < i \leq N$. Because $\left[(\omega-a)^N\right]^{(k)} = P_k^N (\omega-a)^{N-k}$, then $\left[(\omega-a)^N\right]^{(k)}\Big|_{\omega=a} = \begin{cases} 0, & 0 < k < N \\ N! & k = N \end{cases}$.

Furthermore, when $1 \leq i < N$, $f^{(i)}(a) = \left(\sum_{k=0}^{i} C_i^k \left[(\omega-a)^N\right]^{(k)} \varphi^{(i-k)}(\omega)\right)\Big|_{\omega=a} = 0$; when $i = N$,

$f^{(N)}(a) = \left(C_N^N \left[(\omega-a)^N\right]^{(N)} \varphi^{(N-N)}(\omega)\right)\Big|_{\omega=a} = N!\varphi(a)$. Hence, $f^{(i)}(a) = \begin{cases} 0, & 0 \leq i < N \\ N!\varphi(a), & i = N \end{cases}$.

Based on the above Lemma, the proof of Theorem 3 is given as follows.

**Proof of Theorem 3** The following steps are practiced:

(1) It is trivial that $f(d_k) = d_k$ and $f(1) = 1$ under $n \geq 1$, $m \geq 0$ and $m_1 \geq 1$. It is worth noting that $f(\omega)$ satisfies the zero-order condition at $\omega=1$.

(2) Suppose $\varphi(\omega) = (\omega - d_k)^n + c_2(\omega - d_k)(1-\omega)^{m_1} + c_1(1-\omega)^{m+1}$. Because

$\varphi(d_k) = c_1(1-d_k)^{m+1} \neq 0$, then $\left(f(\omega) - d_k\right)^{(i)}\Big|_{\omega=d_k} = \begin{cases} 0, & 0 \leq i \leq n \\ \frac{(n+1)!}{c_1(1-d_k)^{m+1}} (\neq 0), & i = n+1 \end{cases}$ according to above

Lemma.

(3) Through symbolic operation, the first-order differentiation of $f(\omega)$ can be derived as:

$$f'(\omega) = \left(f(\omega) - d_k\right) \left( \frac{n+1}{\omega - d_k} - \frac{n(\omega-d_k)^{n-1} - c_2(\omega-d_k)^{n_1-1}(1-\omega)^{m_1-1} \times \left(\left(n_1(1-\omega) + m_1(\omega-d_k)\right) - c_1(m+1)(1-\omega)^m\right)}{(\omega-d_k)^n + c_2(\omega-d_k)^{n_1}(1-\omega)^{m_1} + c_1(1-\omega)^{m+1}} \right).$$

Consider the situation where $m \geq 1$. Then if $m_1=1$ or/and $m=0$,

$f'(1) = \begin{cases} 1 + c_2(1-d_k)^{1+n_1-n}, & m_1 = 1, m \geq 1 \\ 1 + c_1(1-d_k)^{1-n}, & m = 0, m_1 \geq 2 \\ 1 + c_2(1-d_k)^{1+n_1-n} + c_1(1-d_k)^{1-n}, & m_1 = 0, m = 1 \end{cases}$ ; if $m_1 \geq 2$ and $m \geq 1$, $f'(1) = n + 1 - n = 1$.

We further derive the second-order derivatives of $f(\omega)$, which can be formulated as:

$$f''(\omega) = \frac{h(\omega)}{\left[(\omega-d_k)^n + c_2(\omega-d_k)^{n_1}(1-\omega)^{m_1} + c_1(1-\omega)^{m+1}\right]^3}$$

where $h(\omega) = \sum_i (1-\omega)^{n_h} R_i$ and $R_i$ denotes (rational) polynomial. We have testified that when $R_i$ is rational polynomial, its denominator has the form as $(\omega - d_k)^{n_R}$ with $n_R$ as certain integer. The exponent $n_h$ takes values within $\{m-1, m, 2m-1, 2m; 2m_1-2, 2m_1-1, 2m_1; m_1, m_1+m-2, m_1+m-1; m_1-2, m_1-1, m_1\}$. It can be seen the minimum of $n_h$ will be either $m-1$ or $m_1-2$, and therefore:

(i) If $m_1-2 \geq m-1$ or $m_1 \geq m+1$, the exponent $n_h \geq m-1$, and $f''(\omega)$ can be re-formulated as

$$f''(\omega) = (1-\omega)^{m-1} \varphi(\omega) \quad \text{where} \quad \varphi(1) \neq 0. \quad \text{So} \quad f''(1) = \begin{cases} \varphi(1), & m=1 \\ 0, & m \geq 2 \end{cases}. \quad \text{Let}$$

$$F(\omega) = (-1)^{m-1} f''(\omega) = (\omega-1)^{m-1} \varphi(\omega), \quad \text{then} \quad F^{(i)}(1) = \begin{cases} 0, & 0 \leq i < m-2 \\ (m-1)!\varphi(1), & i = m-1 \end{cases} \quad \text{if} \quad m \geq 2$$

according to the previous Lemma and $F(1) = \varphi(1)$ if $m=1$. One can find by derivation that if $m_1=m+1$, $\varphi(1) = -(1-d_k)^{2n+1} c_1(m+m^2) + (1-d_k)^{2n+1+n_1} c_2(m_1 - m_1^2)$, and if $m_1 > m+1$, $\varphi(1) = -(1-d_k)^{2n+1} c_1(m+m^2)$. Therefore if $m \geq 1$ and $c_1, c_2$ have the same sign, $\varphi(1) \neq 0$, and so

$$f^{(i)}(1) = \begin{cases} (-1)^{m-1} F^{(i-2)}(1) = 0, & 2 \leq i \leq m \\ (-1)^{m-1} F^{(m-1)}(1) \neq 0, & i = m \end{cases} \quad \text{if} \quad m \geq 2 \quad \text{and} \quad f''(1) = (-1)^{m-1} F(1) = \varphi(1) \quad \text{if} \quad m=1. \text{ In}$$

summary, if $c_1$ and $c_2$ have the same sign, $f^{(i)}(1) = \begin{cases} 1 & i=1 \\ 0 & 2 \leq i \leq m \\ \neq 0 & i = m+1 \end{cases}$ if $m \geq 2$ and

$f''(1) = (-1)^{m-1} F(1) = \varphi(1)$ if $m=1$.

(ii) If $m_1-2 < m-1$ under $m_1-3 \geq 0$, or $m+1 > m_1 \geq 3$, the minimum of $n_h$ will be $m_1-2$. Following similar procedures as above, one can find: $\varphi(1) = (1-d_k)^{2n+1+n_1} c_2(m_1 - m_1^2)$. Therefore if $m, m_1 \geq 3$,

$f^{(i)}(1) = \begin{cases} 1 & i=1 \\ 0 & 2 \leq i \leq m_1-1 \\ \neq 0 & i = m_1 \end{cases}$. When $m_1=2$, $n_h=0$, and then $f^{(i)}(1) = \begin{cases} 1 & i=1 \\ \neq 0 & i=2 \end{cases}$.

In short, $f(\omega)$ satisfies $C_{n,\min(m,m_1-1)}$ condition under $d_k \leq \omega \leq 1$, $n \geq 1$, $m \geq 0$ and $m_1 \geq 1$.